\documentclass{amsart}
\usepackage{amsmath,amsthm}
\usepackage{amsfonts,amssymb}
\usepackage[all]{xy}

\usepackage{enumerate}

\newtheorem{theorem}{Theorem}[section]

\newtheorem{lemma}[theorem]{Lemma}
\newtheorem{proposition}[theorem]{Proposition}

\newtheorem{remark}[theorem]{Remark}

\theoremstyle{remark}

\newcommand{\mR}{\mathbb{R}}
\newcommand{\mC}{\mathbb{C}}
\newcommand{\mN}{\mathbb{N}}

\newcommand{\mS}{\mathbb{S}}

\newcommand{\cM}{\mathcal{M}}

\newcommand{\cF}{\mathcal{F}}

\newcommand{\cC}{\mathcal{C}}

\newcommand{\cS}{\mathcal{S}}

\newcommand{\ux}{\underline{x}}

\newcommand{\uy}{\underline{y}}

\newcommand{\upx}{\partial_{\underline{x}}}
\newcommand{\upy}{\partial_{\underline{y}}}

 \def\a{{\alpha}}

 \def\l{{\lambda}}

 \def\la{{\langle}}
 \def\ra{{\rangle}}

 \def\RR{{\mathbb R}}

\newcommand{\wt}{\widetilde}

\begin{document}
 
\title{The class of Clifford-Fourier transforms}

\author{Hendrik De Bie}
\address{Department of Mathematical Analysis\\
Ghent University\\ Krijgslaan 281, 9000 Gent\\ Belgium.}
\email{Hendrik.DeBie@UGent.be}

\author{Nele De Schepper}
\address{Department of Mathematical Analysis\\
Ghent University\\ Krijgslaan 281, 9000 Gent\\ Belgium.}
\email{nds@cage.ugent.be}

\author{Frank Sommen}
\address{Department of Mathematical Analysis\\
Ghent University\\ Krijgslaan 281, 9000 Gent\\ Belgium.}
\email{fs@cage.ugent.be}

\date{\today}
\thanks{H. De Bie is a Postdoctoral Fellow of the Research Foundation - Flanders (FWO)}

\begin{abstract}
Recently, there has been an increasing interest in the study of hypercomplex signals and their Fourier transforms. This paper aims to study such integral transforms from general principles, using 4 different yet equivalent definitions of the classical Fourier transform. This is applied to the so-called Clifford-Fourier transform (see [F. Brackx et al., The Clifford-Fourier transform. {\em J. Fourier Anal. Appl.} {\bf 11} (2005), 669--681]). The integral kernel of this transform is a particular solution of a system of PDEs in a Clifford algebra, but is, contrary to the classical Fourier transform, not the unique solution. Here we determine an entire class of solutions of this system of PDEs, under certain constraints. For each solution, series expressions in terms of Gegenbauer polynomials and Bessel functions are obtained. This allows to compute explicitly the eigenvalues of the associated integral transforms. In the even-dimensional case, this also yields the inverse transform for each of the solutions. Finally, several properties of the entire class of solutions are proven.
\end{abstract}

\maketitle

%\tableofcontents

%%%%%%%%%%%%%%%%%%%%%%%%%%%%%%%%%%%%%%%%%%%%%%%%%%%%%%%%%%%%%%%%%%%%%%%%%%%%%%%%%%%%%%%%%%%%%%%%%%%%%%%%%%%%%%%%%%%%
\section{Introduction}

The last decades, there has been an increasing interest in the theory of hypercomplex signals (i.e. functions taking values in a Clifford algebra) and the possibility of defining and using Fourier transforms that interact with the Clifford algebra structure. This has been investigated from a practical engineering point of view (see e.g. \cite{MR1875365, Ebl1, Ebl2, MR2460142, Fel}) but also from a purely mathematical point of view (see e.g. \cite{MR1887633, MR2200115, MR1308706, MR2226529}) using the function theory of Clifford analysis established in the books \cite{MR697564, MR1169463}. For more references, we refer the reader to the reviewpaper \cite{AIEP}.
Also in applications, there is an increasing interest in having available a good hypercomplex Fourier transform (e.g. in GIS research, see \cite{GIS}).

There are several drawbacks to most of the kernels proposed so far in the literature.  First, several authors work only in low dimensions (dimension 3 or 4, enabling them to use quaternions instead of a full Clifford algebra) which is usually because they have a specific application in mind in these dimensions. Second, and more importantly, most authors use ad hoc formulations for the kernel function of their transforms: they propose very specific kernels, where e.g. the complex unit $I$ is replaced by a generator of the Clifford algebra. Once the kernel is defined, they study in detail all the properties of the related transform. From our perspective, one should work the other way round, namely start from a list of properties or general mathematical principles one wants the transform to satisfy, and then determine all kernels that satisfy these properties.

For that reason, the main aim of this paper is twofold. First of all, we want to use general ideas on Fourier transforms borrowed from other fields of mathematics (in casu the theory of Dunkl operators (see \cite{MR1827871}) and double affine Hecke algebras, the theory of minimal representations) to give a more structural approach to the study of hypercomplex Fourier kernels. We do this by formulating 4 different possible definitions of the classical Fourier transform, and by generalizing these definitions to the Clifford analysis context. 

Secondly, we want to apply these ideas to the so-called Clifford-Fourier transform (introduced in \cite{MR2190678}). This transform was already based on a Lie algebraic approach to the classical Fourier transform, although until recently (see \cite{DBXu}) its kernel was not known in closed form. However, studying this transform using the 4 different definitions mentioned above provides much more insight in this specific transform, and allows us to expand it to a whole class of transforms, all of which will satisfy similar properties (see section \ref{Properties}).

Let us now first give 4 different definitions of the classical Fourier transform, after which we discuss where they appear in the literature (in different fields of mathematics) and what their implications are.

The classical Fourier transform in $\mR^{m}$ can be defined in many ways. In its most basic formulation, it is given by the integral transform
\[
\textbf{F1}\quad  \cF \lbrack f \rbrack (\uy) = \frac{1}{(2 \pi)^{m/2}} \ \int_{\mathbb{R}^m} e^{-I \la \ux,\uy \ra} \ f(\ux) \ dV(\ux)
\]
with $I$ the complex unit, $\la \ux,\uy \ra$ the standard inner product and $dV(\ux)$ the Lebesgue measure on $\mathbb{R}^m$.
Alternatively, one can rewrite the transform as
\[
\textbf{F2}\quad \cF \lbrack f \rbrack (\uy) = \frac{1}{(2 \pi)^{m/2}} \ \int_{\mathbb{R}^m} K(\ux,\uy) \ f(\ux) \ dV(\ux)
\]
where $K(\ux,\uy)$ is, up to a multiplicative constant, the unique solution of the system of PDEs
\[
\partial_{y_{j}} K(\ux,\uy) = - I x_{j} K(\ux,\uy), \quad j =1, \ldots, m.
\]
Yet another formulation is given by
\[
\textbf{F3}\quad \cF = e^{ \frac{I \pi m}{4}} e^{\frac{I \pi}{4}(\Delta - |\ux|^{2})}
\]
with $\Delta$ the Laplacian in $\mR^{m}$. This expression connects the Fourier transform with the Lie algebra $\mathfrak{sl}_{2}$ generated by $\Delta$ and $|\ux|^{2}$ and with the theory of the quantum harmonic oscillator.
Finally, the kernel can also be expressed as an infinite series in terms of special functions as (see \cite[Section 11.5]{MR0010746})
\[
\textbf{F4}\quad K(\ux, \uy) = 2^{\lambda} \Gamma(\lambda)\sum_{k=0}^{\infty}(k+ \lambda) (-I)^{k} (|\ux||\uy|)^{-\lambda} J_{k+ \lambda}(|\ux||\uy|) \; C_{k}^{\lambda}(\langle \underline{\xi},\underline{\eta} \rangle),
\]
where $\underline{\xi}= \ux/|\ux|$, $\underline{\eta} = \uy/|\uy|$ and $\lambda =(m-2)/2$. Here, $J_{\nu}$ is the Bessel function and $C_{k}^{\l}$ the Gegenbauer polynomial.

Each formulation has its specific advantages and uses. The classical formulation $\textbf{F1}$ allows to immediately compute a bound of the kernel and is hence ideal to study the transform on $L_{1}$ spaces or more general function spaces. 

Formulation $\textbf{F2}$ yields the calculus properties of the transform, and allows to generalize the transform to e.g. the so-called Dunkl transform (see \cite{deJ}). This formulation (defining the kernel as the solution of an eigenvalue problem) is also frequently used in the context of double affine Hecke algebras (see e.g. \cite{C, J}), an algebraic generalization of Dunkl operators. 

Formulation $\textbf{F3}$ emphasizes the structural (Lie algebraic) properties of the Fourier transform and also allows to compute its eigenfunctions and spectrum. This formulation stems from representation theory (see \cite{MR0983366, MR0974332}) and has been used in recent work on minimal representations (see \cite{MR2134314, MR2401813} and further generalizations in \cite{Orsted2}).

Finally, $\textbf{F4}$ connects the Fourier transform with the theory of special functions, and is the ideal formulation to obtain e.g. the Bochner identities (which are a special case of the subsequent Proposition \ref{Bochner}). Similar series representations have also been used in the context of Dunkl operators and have applications in the study of generalized translation operators (see e.g. \cite{MR1973996, DBXu}).

In \cite{MR2190678}, $\textbf{F3}$ was adapted to the case of functions taking values in the Clifford algebra $\cC l_{0,m}$ to define a couple of Fourier transforms in Clifford analysis by
\[
\cF_{\pm} =  e^{ \frac{I \pi m}{4}} e^{\mp \frac{I \pi}{2}\Gamma  }e^{\frac{I \pi}{4}(\Delta - |\ux|^{2})} = e^{ \frac{I \pi m}{4}} e^{\frac{I \pi}{4}(\Delta - |\ux|^{2} \mp 2\Gamma)}
\]
with $2\Gamma= (\upx \ux - \ux \upx)+m$. Here, $\upx$ is the Dirac operator and $\ux$ the vector variable. The exponential now contains the generators of the Lie superalgebra $\mathfrak{osp}(1|2)$. For several years, the problem remained open to write this Clifford-Fourier transform as an integral transform and to determine explicitly its kernel. A breakthrough was obtained in \cite{DBXu}, where the kernel was determined in all even dimensions, and the problem for odd dimensions was reduced to dimension 3 (where an integral representation of the kernel was obtained).

As a by-product, it was also obtained that the kernel $K_+(\ux,\uy)$ of the integral transform $\cF_{+}$ satisfies a system of PDEs, namely
\begin{align}
\label{CFstelselintro}
\begin{split}
\partial_{\uy} \lbrack K_+(\ux,\uy) \rbrack = (-I)^m \ \bigl( K_+(\ux,-\uy) \bigr)^c \ \ux\\
\lbrack K_+(\ux,\uy) \rbrack \partial_{\ux}  =  (-I)^m \ \uy \ \bigl( K_+(\ux,-\uy) \bigr)^c,
\end{split}
\end{align}
where $c$ denotes the complex conjugation.

The main aim of this paper is to study this system of PDEs. We will show that, contrary to the classical formulation $\textbf{F2}$, this system does not have a unique solution, but instead $m-1$ linearly independent solutions $K_{+,m}^{i}$ (when we restrict to a special subclass of solutions satisfying nice symmetries). Each of these solutions gives rise to an associated integral transform
\[
\mathcal{F}^{i}_{+,m} \lbrack f(\ux) \rbrack (\uy) = \frac{1}{(2 \pi)^{m/2}} \  \int_{\mathbb{R}^m} K^{i}_{+,m}(\ux,\uy) \ f(\ux) \ dV(\ux)
\]
and we study each of these transforms in-depth. In particular, we determine series representations of the form $\textbf{F4}$ for all relevant solutions of (\ref{CFstelselintro}). This in turn allows us to obtain the spectrum for the associated integral transforms and allows us to prove the surprising fact that in case of $m$ even
\[
\mathcal{F}^{i}_{+,m} \mathcal{F}^{m-2-i}_{+,m} = \mbox{id}.
\]
In other words, for $m$ even we find a complete class of integral transforms, where the inverse of each element is again an element of the class.
We also obtain bounds on the kernels $K_{+,m}^{i}$, which allows us to define the broadest function space on which the associated transform is defined (compare with $\textbf{F1}$). Finally, we prove several important properties for all the kernels obtained.

%%%%%%%%%%%%%

One of the strengths of our results is that for the kernels obtained in this paper, we obtain always both the $\textbf{F4}$ and $\textbf{F1}$ formulation. This is much more than in, say, the Dunkl case, where for almost all finite reflection groups the formulation $\textbf{F1}$ is missing and one has to use different and complicated techniques to prove e.g. boundedness of the transform.

%%%%%%%%%%%%%

The paper is organized as follows. In section \ref{Preliminaries} we repeat basic notions of Clifford algebras and related differential operators. We give the explicit expression of the kernel of the Clifford-Fourier transform and of the Fourier-Bessel transform. In section \ref{Series} we prove some general statements for kernels expressed as series of products of Gegenbauer polynomials and Bessel functions. In section \ref{Even} we study the Clifford-Fourier system (\ref{CFstelselintro}) in even dimension. We determine an interesting class of solutions, find recursion relations between these solutions and obtain series expansions. We also determine the eigenvalues for each solution. In section \ref{Odd} we treat the case of odd dimension. We omit most proofs in this section, because they are similar as in the even dimensional case. Nevertheless, this case has to be considered separately, because the solutions will now be complex instead of real. Finally, in section \ref{Properties}, we collect some important properties of the new class of Clifford-Fourier transforms and prove the important fact that in the even dimensional case also the inverse of each transform is again an element of the same class.

%%%%%%%%%%%%%%%%%%%%%%%%%%%%%%%%%%%%%%%%%%%%%%%%%%%%%%%%%%%%%%%%%%%%%%%%%%%%%%%%%%%%%%%%%%%%%%%%%%%%%%%%%%%%%%%%%%
%%%%%%%%%%%%%%%%%%%%%%%%%%%%%%%%%%%%%%%%%%%%%%%%%%%%%%%%%%%%%%%%%%%%%%%%%%%%%%%%%%%%%%%%%%%%%%%%%%%%%%%%%%%%%%%%%%%%%
\section{Preliminaries}
\setcounter{equation}{0}
\label{Preliminaries}

\subsection{Clifford analysis and special functions}
Clifford analysis (see e.g. \cite{MR1169463}) is a theory that offers a natural generalization of complex analysis to higher dimensions. To  $\mR^{m}$, the Euclidean space in $m$ dimensions, we first associate the Clifford algebra $\cC l_{0,m}$, generated by the canonical basis $e_{i}$, $i= 1, \ldots, m$. These generators satisfy the multiplication rules $e_{i} e_{j} + e_{j} e_{i} = - 2 \delta_{ij}$.

The Clifford algebra $\cC l_{0,m}$ can be decomposed as $\cC l_{0,m} = \oplus_{k=0}^{m} \cC l_{0,m}^{k}$
with $\cC l_{0,m}^{k}$ the space of $k$-vectors defined by
\[
\cC l_{0,m}^{k} = \mbox{span} \{ e_{i_{1} \ldots i_{k}} = e_{i_{1}} \ldots e_{i_{k}}, i_{1} < \ldots < i_{k} \}.
\]
More precisely, we have that the space of $1$-vectors is given by $\cC l_{0,m}^{1} = \mbox{span} \{ e_{i}, i = 1, \ldots, m \}$ 
and it is obvious that this space is isomorphic with $\mR^{m}$. The space of so-called bivectors is given explicitly by $\cC l_{0,m}^{2} = \mbox{span} \{ e_{ij}=e_{i} e_{j}, i < j \}$.

We identify the point $(x_{1}, \ldots, x_{m})$ in $\mR^{m}$ with the vector variable $\ux$ given by $\ux = \sum_{j=1}^{m} x_{j} e_{j}$. The Clifford product of two vectors splits into a scalar part and a bivector part:
\[
\ux \uy = \ux . \uy + \ux \wedge \uy,
\]
with
\[
\ux . \uy = - \langle \ux, \uy \rangle = -\sum_{j=1}^{m} x_{j} y_{j} = \frac{1}{2} (\ux \uy + \uy  \ux)
\]
and
\[
\ux \wedge \uy = \sum_{j<k} e_{jk} (x_{j} y_{k} - x_{k}y_{j}) = \frac{1}{2} (\ux \uy - \uy  \ux).
\]
It is interesting to note that the square of a vector variable $\ux$ is scalar-valued and equals the norm squared up to a minus sign: $
\ux^{2} = - \langle \ux, \ux \rangle = - |\ux|^{2}$.
Similarly, we introduce a first order vector differential operator by
\[
\upx = \sum_{j=1}^{m} \partial_{x_{j}} e_{j}.
\]
This operator is the so-called Dirac operator. Its square equals, up to a minus sign, the Laplace operator in $\mR^{m}$: $\upx^{2} = - \Delta$. A function $f$ defined in some open domain $\Omega \subset \mR^{m}$ with values in the Clifford algebra $\cC l_{0,m}$ is called monogenic if $\upx ( f ) = 0$.

Another important operator in Clifford analysis is the so-called Gamma operator, defined by
\[
\Gamma_{\ux} =  - \ux \wedge \upx = - \sum_{j<k} e_{jk} (x_{j} \partial_{x_{k}} - x_{k}\partial_{x_{j}}).
\]
This operator is bivector-valued.

A basis $\lbrace \psi_{j,k,\ell} \rbrace$ for the space $\mathcal{S}(\mathbb{R}^m) \otimes \cC l_{0,m}$, where $\mathcal{S}(\mathbb{R}^m)$ denotes the Schwartz space, is given by (see \cite{MR926831})
\begin{align} \label{basis}
\begin{split}
\psi_{2j,k,\ell}(\ux) &:= L_{j}^{\frac{m}{2}+k-1}(|\ux|^{2}) \ M_{k}^{(\ell)}(\ux) \ e^{-|\ux|^{2}/2},\\
\psi_{2j+1,k,\ell}(\ux) &:= L_{j}^{\frac{m}{2}+k}(|\ux|^{2}) \ \ux \ M_{k}^{(\ell)}(\ux) \ e^{-|\ux|^{2}/2},
\end{split}
\end{align}
where $j,k \in \mathbb{N}$, $L_j^{\alpha}$ are the Laguerre polynomials and $ \{M_k^{(\ell)}\}$, $(\ell=1,2, \ldots , \mathrm{dim}(\mathcal{M}_k))$ is a basis for the space $\mathcal{M}_k$. $\mathcal{M}_k$ is the space of spherical monogenics of degree $k$, i.e. homogeneous polynomial null-solutions of the Dirac operator of degree $k$.

In the sequel we will frequently need the following well-known properties of Gegenbauer polynomials (see e.g. \cite{Sz}):
\begin{equation}\label{Geg1}
\frac{\lambda + n}{\lambda} \ C_n^{\lambda}(w) = C_n^{\lambda + 1}(w) - C_{n-2}^{\lambda+1}(w)
\end{equation}
and
\begin{equation}\label{Geg2}
w \ C_{n-1}^{\lambda+1}(w) = \frac{n}{2(n+\lambda)} \ C_n^{\lambda+1}(w) + \frac{n+ 2 \lambda}{2(n+\lambda)} \ C_{n-2}^{\lambda+1}(w),
\end{equation}
as well as the Bessel function identity
\begin{equation}\label{Besselid}
J_{\nu}(z) = \frac{z}{2 \nu} \ \left(J_{\nu+1}(z) + J_{\nu-1}(z) \right).
\end{equation}
%%%%%%%%%%%%%%%%%%%%%%%%%%%%%%%%%%%%%%%%%%%%%%%%%%%%%%%%%%%%%%%%%%%%%%%%%%%%%%%%%%%%%%%%%%
\subsection{The Clifford-Fourier transform}
Several attempts have been made to introduce a generalization of the classical Fourier transform \textbf{F1} to the setting of Clifford analysis (see the introduction and \cite{AIEP} for a review). We will concentrate on the so-called Clifford-Fourier transform introduced in \cite{MR2190678} by an operator exponential, similar as the \textbf{F3} representation of the classical Fourier transform:
\[
\cF_{\pm} =  e^{ \frac{I \pi m}{4}} e^{\mp \frac{I \pi}{2}\Gamma  }e^{\frac{I \pi}{4}(\Delta - |\ux|^{2})}. 
\]

This Fourier type transform can equivalently be written as an integral transform
\begin{displaymath}
\cF_{\pm} \lbrack f \rbrack (\uy) = \frac{1}{(2 \pi)^{m/2}} \ \int_{\mathbb{R}^m} K_{\pm}(\ux,\uy) \ f(\ux) \ dV(\ux),
\end{displaymath}
where the kernel function $K_{\pm}(\ux,\uy)$ is given by the operator exponential $e^{\mp I \frac{\pi}{2} \Gamma}$ acting on the classical Fourier kernel, i.e.
\begin{equation}
\label{Kexpop}
K_{\pm}(\ux,\uy) = e^{\mp I \frac{\pi}{2} \Gamma_{\uy}} \left( e^{-I \la\ux,\uy\ra} \right) .
\end{equation}

%\begin{remark}
%As the kernel is not symmetric, i.e. $K_{\pm}(\ux,\uy) \neq K_{\pm}(\uy,\ux)$, we will always assume in this paper that we integrate over the first variable in the kernel.
%\end{remark}

Similar to the classical case, the Clifford-Fourier transform satisfies some calculus rules, which translates to the following system of equations satisfied by the kernel:
\begin{align} \label{C-F system}
\begin{split}
\partial_{\uy} \lbrack K_{\mp}(\ux,\uy) \rbrack = \mp \ (\pm I)^m \ K_{\pm}(\ux,\uy) \ \ux\\
\lbrack K_{\pm}(\ux,\uy) \rbrack \partial_{\ux} = \pm (\mp I)^m \ \uy \ K_{\mp}(\ux,\uy),
\end{split}
\end{align}
where 
\[
\lbrack K_{\pm}(\ux,\uy) \rbrack \partial_{\ux}  = \sum_{i=1}^{m} \left(\partial_{x_{i}}K_{\pm}(\ux,\uy) \right) e_{i}
\]
denotes the action of the Dirac operator on the right. The system of PDEs (\ref{C-F system}) should be compared with the formulation \textbf{F2} of the classical Fourier transform.

Explicit computation of (\ref{Kexpop}) is a hard problem. Until recently, the Clifford-Fourier kernel was known explicitly only in the case $m=2$ (see \cite{MR2283868}); for higher even dimensions, a complicated iterative procedure for constructing the kernel was given in \cite{JFAA-Fourier-Bessel}, which could only be used practically in low dimensions. A breakthrough was obtained in \cite{DBXu}. In this paper it is found that for $m$ even the kernel can be expressed as follows in terms of a finite sum of Bessel functions:
\begin{equation}\label{kernel C-F}
K_{+}(\ux,\uy) =  \left( \frac{\pi}{2} \right)^{1/2} \ \left( A(s,t) + B(s,t) + (\ux \wedge \uy) \ C(s,t) \right)
\end{equation}
with
\begin{align}
\label{CFeven}
\begin{split}
A(s,t) &=  \sum_{\ell=0}^{\lfloor \frac{m}{4}-\frac{3}{4} \rfloor} s^{m/2-2-2 \ell} \ \frac{1}{2^{\ell} \ell!} \frac{\Gamma \left( \frac{m}{2} \right)}{\Gamma \left( \frac{m}{2}-2 \ell -1 \right)} \ \widetilde{J}_{(m-2\ell-3)/2}(t)\\
B(s,t) &=   \sum_{\ell=0}^{\lfloor \frac{m}{4}-\frac{1}{2} \rfloor} s^{m/2-1-2 \ell} \ \frac{1}{2^{\ell} \ell!} \frac{\Gamma \left( \frac{m}{2} \right)}{\Gamma \left( \frac{m}{2}-2 \ell  \right)} \ \widetilde{J}_{(m-2\ell-3)/2}(t)\\
C(s,t) &=  - \sum_{\ell=0}^{\lfloor \frac{m}{4}-\frac{1}{2} \rfloor} s^{m/2-1-2 \ell} \ \frac{1}{2^{\ell} \ell!} \frac{\Gamma \left( \frac{m}{2} \right)}{\Gamma \left( \frac{m}{2}-2 \ell \right)} \ \widetilde{J}_{(m-2\ell-1)/2}(t).
\end{split}
\end{align}
Here $\lfloor \ell \rfloor$ denotes the largest $n \in \mathbb{N}$ which satisfies $n \leq \ell$ and the notations $s= \langle \ux,\uy \rangle$, $t=| \ux \wedge \uy|= \sqrt{|\ux|^{2} |\uy|^{2}-\langle \ux, \uy \rangle^2}$ and $\widetilde{J}_{\alpha}(t) = t^{-\alpha} J_{\alpha}(t)$ are used. Moreover, it is shown that
\begin{equation}\label{rel C-F kernel}
K_+(\ux,\uy) = \bigl( K_{-}(\ux,-\uy) \bigr)^c
\end{equation}
holds and also that in the case $m$ even, $K_{-}(\ux,\uy)$ is real-valued, hence in this case the complex conjugation in the above relation can be omitted. Note that the Clifford-Fourier kernel is parabivector-valued, i.e. it takes the form of a scalar plus a bivector.

For $m$ odd, the question of determining the kernel explicitly was reduced to the case of $m=3$. There, a more or less complicated integral expression of the kernel was obtained (see \cite[Lemma 4.5]{DBXu}). A simple expression as in formula (\ref{CFeven}) is not known in this case.

Finally, let us mention the action of the Clifford-Fourier transform on the basis elements (\ref{basis}) (see \cite{MR2190678}):
\begin{align} \label{eigenvalue C-F}
\begin{split}
\cF_{\pm}\lbrack \psi_{2p,k,\ell}\rbrack(\uy) &= (-1)^{p+k} \ (\mp 1)^k \ \psi_{2p,k,\ell}(\uy)\\
\cF_{\pm}\lbrack \psi_{2p+1,k,\ell}\rbrack(\uy) &= I^m \ (-1)^{p+1} \ (\mp1)^{k+m-1} \ \psi_{2p+1,k,\ell}(\uy).
\end{split}
\end{align}
%%%%%%%%%%%%%%%%%%%%%%%%%%%%%%%%%%%%%%%%%%%%%%%%%%%%%%%%%%%%%%%%%%%%%%%%%%%%%%%%%%%%%%%%%%%5
\subsection{The Fourier-Bessel transform}
In \cite{FourierBessel} another new integral transform within the Clifford analysis setting was devised, the so-called Fourier-Bessel transform given by
\begin{displaymath}
\cF^{\mathrm{Bessel}} \lbrack f \rbrack (\uy) = \frac{1}{(2 \pi)^{m/2}} \ \int_{\mathbb{R}^m} K^{\mathrm{Bessel}}(\ux,\uy) \ f(\ux) \ dV(\ux).
\end{displaymath}
Its integral kernel takes the form
\begin{equation}\label{kernel F-B}
K^{\mathrm{Bessel}}(\ux,\uy) = \sqrt{\frac{\pi}{2}} \ \left( (-1)^{m/2} \ \widetilde{J}_{(m-3)/2}(t) + (\ux \wedge \uy) \ \widetilde{J}_{(m-1)/2}(t) \right),
\end{equation}
where we use a different normalization as in \cite{FourierBessel}. 

Note that similar to the Clifford-Fourier kernel, it is parabivector-valued. Moreover, the basis elements (\ref{basis}) are also eigenfunctions of this transform. The eigenvalues are however quite a bit more complicated in this case. More precisely, we have for $k$ even
\begin{align} \label{eigenvalue F-B1}
\begin{split}
\cF^{\mathrm{Bessel}}\lbrack \psi_{2p,k,\ell}\rbrack(\uy) &= (-1)^{m/2} \ (-1)^p \ \frac{(k-1)!!}{(k+m-3)!!} \ \psi_{2p,k,\ell}(\uy)\\
\cF^{\mathrm{Bessel}}\lbrack \psi_{2p+1,k,\ell}\rbrack(\uy) &= (-1)^p \ \frac{(k-1)!!}{(k+m-3)!!} \ \psi_{2p+1,k,\ell}(\uy),
\end{split}
\end{align}
while for $k$ odd
\begin{align} \label{eigenvalue F-B2}
\begin{split}
\cF^{\mathrm{Bessel}}\lbrack \psi_{2p,k,\ell}\rbrack(\uy) &= (-1)^{p+1} \ \frac{k!!}{(k+m-2)!!} \  \psi_{2p,k,\ell}(\uy)\\
\cF^{\mathrm{Bessel}}\lbrack \psi_{2p+1,k,\ell}\rbrack(\uy) &= (-1)^{m/2} \ (-1)^p \ \frac{k!!}{(k+m-2)!!} \ \psi_{2p+1,k,\ell}(\uy).
\end{split}
\end{align}
For $u$ odd, $u!!$ denotes the product: $u!!=u(u-2)(u-4) \ldots 5 \ 3 \ 1$, while for $u$ even, $u!!$ stands for the product: $u!!=u(u-2) \ldots 6 \ 4 \ 2$.
%%%%%%%%%%%%%%%%%%%%%%%%%%%%%%%%%%%%%%%%%%%%%%%%%%%%%%%%%%%%%%%%%%%%%%%%%%%%%%%%%%%%%%%%%%%%%%%%%%%%%%%%%%%%
%%%%%%%%%%%%%%%%%%%%%%%%%%%%%%%%%%%%%%%%%%%%%%%%%%%%%%%%%%%%%%%%%%%%%%%%%%%%%%%%%%%%%%%%%%%%%%%%%%%%%%%%%%%%
\section{Series approach}
\setcounter{equation}{0}
\label{Series}

In this section we consider a general kernel of the following form
\begin{equation}\label{structure kernel}
K_{-}(\ux,\uy) = A(w,z) + (\ux \wedge \uy) \  B(w,z)
\end{equation}
with
\begin{align}
\label{seriesconv}
\begin{split}
A(w,z) &= \sum_{k=0}^{+\infty} \alpha_{k}z^{-\l}J_{k+\l}(z) C^{\l}_{k}(w)\\
B(w,z) &=  \sum_{k=1}^{+\infty} \beta_{k}z^{-\l-1}J_{k+\l}(z) C^{\l+1}_{k-1}(w)
\end{split}
\end{align}
and $\alpha_{k}, \beta_{k} \in \mC$. Here, we have introduced the variables $z= |\ux| |\uy|$, $w=\langle \underline{\xi},\underline{\eta} \rangle$ ($\ux = |\ux| \underline{\xi}$, $\uy = |\uy| \underline{\eta}$, $\underline{\xi},\underline{\eta} \in \mS^{m-1}$) and use the notation $\lambda=(m-2)/2$. The kernel $K_{+}(\ux,\uy)$ is then obtained by the formula $K_{+}(\ux,\uy) = \bigl( K_{-}(\ux,-\uy) \bigr)^c$.

Note that the convergence of the series in (\ref{seriesconv}) is never a problem for the coefficients $\a_{k}$ and $\beta_{k}$ we will consider. Indeed, we can e.g. estimate
\begin{align*}
\left| \sum_{k=0}^{+\infty} \alpha_{k}z^{-\l}J_{k+\l}(z) C^{\l}_{k}(w)\right| &\leq 2^{-\l} \sum_{k=0}^{+\infty} |\alpha_{k}| \left|\left(\frac{z}{2}\right)^{-\l-k}J_{k+\l}(z) \right| \left(\frac{z}{2}\right)^{k} |C^{\l}_{k}(w)| \\
& \leq 2^{-\l} \l B(\l)\sum_{k=0}^{+\infty} |\alpha_{k}| \frac{1}{\Gamma(k+\l+1)} \left(\frac{z}{2}\right)^{k} k^{2\l-1} 
\end{align*}
where we used the estimate 
\[
\left|\left(\frac{z}{2}\right)^{-\l-k}J_{k+\l}(z) \right| \leq \frac{1}{\Gamma(k+\l+1)}
\]
which follows immediately from the integral representation of the Bessel function (see \cite[(1.71.6)]{Sz})
and the fact that there exists a constant $B(\l)$ such that
\begin{equation*}
\sup_{-1\leq w \leq 1} \left|\frac{1}{\l}C_k^{\l}(w) \right| \leq B(\l) k^{2\l-1}, \quad \forall k \in \mN,
\end{equation*}
see \cite[Lemma 4.9]{Orsted2}. We conclude that if $\alpha_{k}$ is a fixed rational function of $k$ (as will be the case in Theorem \ref{SeriesEven} and \ref{SeriesOdd}), then the series converges absolutely and uniformly on compacta because of the ratio test.

We define the following two integral transforms
\begin{align*}
\cF_{\pm}\lbrack f \rbrack (\uy) &= \frac{\Gamma \left( \frac{m}{2} \right)}{2 \pi^{m/2}} \int_{\mR^{m}}K_{\pm}(\ux,\uy)\ f(\ux) \ dV(\ux).
\end{align*}
Now we calculate the action of these transforms on the basis (\ref{basis}) of $\cS(\mR^{m}) \otimes \cC l_{0,m}$. We start with the following auxiliary result, which is a generalization of the Bochner formulas for the classical Fourier transform.

\begin{proposition}
\label{Bochner}
Let $M_{k} \in \cM_{k}$ be a spherical monogenic of degree $k$. Let $f(\ux)= f_0(|\ux|)$ be a real-valued radial function in 
$\cS(\RR^m)$. Further, put $\underline{\xi}= \ux/|\ux|$, $\underline{\eta} = \uy/|\uy|$ and $r = |\ux|$. Then one has
\begin{eqnarray*} 
\cF_{-} \left\lbrack f(\ux)M_{k}(\ux) \right\rbrack (\uy) &=& \left( \frac{\l}{\l+k} \alpha_{k} - \frac{k}{2(k+ \l)} \beta_k \right)  M_{k}(\underline{\eta})\\
&&\times \int_{0}^{+\infty} r^{m+k-1}f_0(r)   z^{-\l} J_{k + \l}(z)  dr
\end{eqnarray*}
and
\begin{eqnarray*}
\cF_{-} \left\lbrack f(\ux) \ux M_{k}(\ux) \right\rbrack (\uy) &=& \left( \frac{\l}{\l+k+1} \alpha_{k+1} + \frac{k+1+2\l}{2(k+1+ \l)} \beta_{k+1} \right)  \underline{\eta} \  M_{k}(\underline{\eta})\\
&&\times \int_{0}^{+\infty} r^{m+k}f_0(r)   z^{-\l} J_{k +1+ \l}(z)  dr
\end{eqnarray*}
with $z= r |\uy|$ and $\l = (m-2)/2$.
\end{proposition}

\begin{proof}
The proof goes along similar lines as the proof of Theorem 6.4 in \cite{DBXu}.
\end{proof}

We then have the following theorem.

\begin{theorem}
\label{eigenvalues}
One has, putting $\beta_{0}=0$,
\begin{align}\label{eigenvalue eq}
\begin{split}
\cF_{-}\lbrack \psi_{2j,k,\ell}\rbrack (\uy) & = \left(\frac{\l}{\l+k} \alpha_{k} - \frac{k}{2(\l+k)} \beta_{k} \right)(-1)^{j}\psi_{2j,k,\ell}(\uy)\\
\cF_{-}\lbrack \psi_{2j+1,k,\ell}\rbrack (\uy) & = \left(\frac{\l}{\l+k+1} \alpha_{k+1} + \frac{k+1+2\l}{2(\l+k+1)} \beta_{k+1} \right)(-1)^{j}\psi_{2j+1,k,\ell}(\uy)
\end{split}
\end{align}
and 
\begin{align*}
\cF_{+}\lbrack \psi_{2j,k,\ell}\rbrack (\uy) & = \left(\frac{\l}{\l+k} \alpha_{k}^c - \frac{k}{2(\l+k)} \beta_{k}^c \right)(-1)^{j+k}\psi_{2j,k,\ell}(\uy)\\
\cF_{+}\lbrack \psi_{2j+1,k,\ell}\rbrack (\uy) & = \left(\frac{\l}{\l+k+1} \alpha_{k+1}^c + \frac{k+1+2\l}{2(\l+k+1)} \beta_{k+1}^c \right)(-1)^{j+k+1}\psi_{2j+1,k,\ell}(\uy).
\end{align*}
\end{theorem}

\begin{proof}
This follows from the explicit expression (\ref{basis}) of the basis and the identity (see e.g. \cite[exercise 21, p. 371]{Sz})
\[
\int_{0}^{+\infty} r^{2\l+1} (rs)^{-\l} J_{k+\l}(rs)\, r^{k} L_{j}^{k+\l}(r^{2}) e^{-r^{2}/2}dr = (-1)^{j}s^{k} L_{j}^{k+\l}(s^{2}) e^{-s^{2}/2}. 
\]
\end{proof}

We are now able to construct the inverse of $\cF_{-}$ on the basis $\{\psi_{j,k,\ell}\}$. The construction is similar for $\cF_{+}$.

\begin{theorem}
The inverse of $\cF_{-}$ on the basis $\{\psi_{j,k,\ell}\}$ is given by
\[
\cF_{-}^{-1}\lbrack f \rbrack (\uy) =  \frac{\Gamma \left( \frac{m}{2} \right) }{2 \pi^{m/2}} \int_{\mR^{m}}\wt{K_{-}(\ux,\uy)} \ f(\ux) \ dV(\ux)\\
\]
with $\wt{K_{-}(\ux,\uy)} = A(w,z) + (\ux \wedge \uy) \  B(w,z)$ given by
\begin{align*}
A(w,z)&= \sum_{k=0}^{+\infty} \frac{1}{N}(\alpha_{k}+ \beta_{k}) z^{-\l}J_{k+\l}(z) C^{\l}_{k}(w)\\
B(w,z) &=  - \sum_{k=1}^{+\infty} \frac{1}{N} \beta_{k}z^{-\l-1}J_{k+\l}(z) C^{\l+1}_{k-1}(w),
\end{align*}
where
\[
N = \left(\frac{\l}{\l+k} \alpha_{k} - \frac{k}{2(\l+k)} \beta_{k} \right)\left(\frac{\l}{\l+k} \alpha_{k} + \frac{k+2\l}{2(\l+k)} \beta_{k} \right).
\]
\end{theorem}

\begin{proof}
Put $\wt{K_{-}(\ux,\uy)} = A(w,z) + (\ux \wedge \uy) \  B(w,z)$ where
\begin{align*}
A(w,z) &= \sum_{k=0}^{+\infty} \gamma_{k} z^{-\l}J_{k+\l}(z) C^{\l}_{k}(w)\\
B(w,z) &=  \sum_{k=1}^{+\infty}\delta_{k} z^{-\l-1}J_{k+\l}(z) C^{\l+1}_{k-1}(w)
\end{align*}
and with $\gamma_{k}, \delta_{k} \in \mC$. We need to have that
\[
\cF_{-}^{-1} \bigl\lbrack \cF_{-} \lbrack f \rbrack \bigr\rbrack = \cF_{-} \bigl\lbrack \cF_{-}^{-1} \lbrack f \rbrack \bigr\rbrack = f.
\]
Using Theorem \ref{eigenvalues} this condition is equivalent with the system of equations ($k=0,1, \ldots$)
\begin{align*}
&\left(\frac{\l}{\l+k} \alpha_{k} - \frac{k}{2(\l+k)} \beta_{k} \right) \left(\frac{\l}{\l+k} \gamma_{k} - \frac{k}{2(\l+k)} \delta_{k} \right) =1\\
&\left(\frac{\l}{\l+k} \alpha_{k} + \frac{k+2\l}{2(\l+k)} \beta_{k} \right)\left(\frac{\l}{\l+k} \gamma_{k} + \frac{k+2\l}{2(\l+k)} \delta_{k} \right)=1.
\end{align*}
Solving this system then yields the statement of the theorem.
\end{proof}

Now our aim is to see what restrictions should be put on the coefficients $\alpha_k$ and $\beta_k$ such that $\cF_{\pm}$ satisfies the Clifford-Fourier system:
\begin{align*}
\cF_{\pm} \left\lbrack \ux \, f \right\rbrack(\uy) & = \mp (\mp I)^{m} \upy \left\lbrack \cF_{\mp} \left\lbrack f \right\rbrack (\uy) \right\rbrack \\
\cF_{\pm} \left\lbrack \upx \lbrack f \rbrack  \right\rbrack (\uy) & = \mp (\mp I)^{m}  \uy \ \cF_{\mp} \left\lbrack f \right\rbrack (\uy)
\end{align*}
and more specifically
\begin{equation}\label{recursionFourier}
\cF_{\pm} \left\lbrack(\upx-\ux) \lbrack f \rbrack \right\rbrack (\uy) = \pm (\mp I)^{m}  ( \upy - \uy) \left\lbrack \cF_{\mp} \left\lbrack f \right\rbrack (\uy) \right\rbrack.
\end{equation}
Now recall that (see \cite{MR926831})
\begin{displaymath}
\psi_{j,k,\ell}(\ux) = \frac{(-1)^j \ 2^{-j}}{\left\lfloor \frac{j}{2} \right\rfloor!} \ (\upx-\ux)^{j} \left\lbrack M_{k}^{(\ell)}(\ux) \  e^{-r^{2}/2} \right\rbrack.
\end{displaymath}
We then have, on the one hand, using Theorem \ref{eigenvalues}
\[
\cF_{+}\left\lbrack \psi_{2j+1,k,\ell}\right\rbrack (\uy) = \left(\frac{\l}{\l+k+1} \alpha_{k+1}^c + \frac{k+1+2\l}{2(\l+k+1)} \beta_{k+1}^c \right) \ (-1)^{j+k+1} \ \psi_{2j+1,k,\ell}(\uy)
\]
and on the other hand, using (\ref{recursionFourier})
\begin{align*}
\cF_{+} \left\lbrack \psi_{2j+1,k,\ell}\right\rbrack (\uy) &= - \frac{1}{2} \ \cF_{+}\left\lbrack (\upx-\ux) \lbrack \psi_{2j,k,\ell}\rbrack \right\rbrack (\uy)\\
&=  -\frac{1}{2} \ (- I)^{m}  \ ( \upy - \uy)  \left\lbrack \cF_{-} \lbrack \psi_{2j,k,\ell} \rbrack (\uy) \right\rbrack\\
&=  (- I)^{m}  \ \left(\frac{\l}{\l+k} \alpha_{k} - \frac{k}{2(\l+k)} \beta_{k} \right) \ (-1)^{j} \ \psi_{2j+1,k,\ell}(\uy).
\end{align*}
This leads to the following condition on $\alpha_{k}$ and $\beta_{k}$:
\[
 \left(\l \alpha_{k+1}^c + \frac{k+1+2\l}{2} \beta_{k+1}^c \right) =  (- I)^{m}(-1)^{k+1} \frac{\l+k+1}{\l+k}  \left(\l \alpha_{k} - \frac{k}{2} \beta_{k} \right).
\]
%%%%%%%%%%%%%%%%%%%%%%%%%%%%%%%%%%%%%%%%%%%%%%%%%%%%%%%%%%%%%%%%%%%%%%%%%%%%%%%%%%%%%%%%%%%%%%%%%%%%%%%%%%%
%%%%%%%%%%%%%%%%%%%%%%%%%%%%%%%%%%%%%%%%%%%%%%%%%%%%%%%%%%%%%%%%%%%%%%%%%%%%%%%%%%%%%%%%%%%%%%%%%%%%%%%%%%
\section{New Clifford-Fourier transforms: the case $m$ even}
\setcounter{equation}{0}
\label{Even}

\subsection{Parabivector-valued solutions of the Clifford-Fourier system}
The aim of this section is to solve the Clifford-Fourier system (\ref{C-F system}) in even dimension:
\begin{align}\label{system}
\begin{split}
\partial_{\uy} \lbrack K^+(\ux,\uy) \rbrack  =  a \ K^-(\ux,\uy) \ \ux\\
\lbrack K^+(\ux,\uy) \rbrack \partial_{\ux}  =  a \ \uy \ K^-(\ux,\uy)
\end{split}
\end{align}
with $ a=(-1)^{m/2}$ and $K^-(\ux,\uy) = K^+(\ux,-\uy)$, see (\ref{rel C-F kernel}).

As the even dimensional Clifford-Fourier transform is real-valued, we look for real-valued solutions $K^+(\ux,\uy)$. Inspired by the expression (\ref{kernel C-F}), we want to determine parabivector-valued solutions of the form:
\begin{align*}
K^+(\ux,\uy) &= f(s,t) + (\ux \wedge \uy) \ g(s,t)\\
K^-(\ux,\uy) &=  f(-s,t) - (\ux \wedge \uy) \ g(-s,t)
\end{align*}
with $s = \langle \ux,\uy \rangle$, $t=|\ux \wedge \uy|$ and $f$ and $g$ real-valued functions.

Taking into account that (see e.g. \cite{FourierBessel})
\begin{displaymath}
\partial_{\uy} \lbrack s \rbrack = \ux \ \ \ , \ \ \ \partial_{\uy}\lbrack t \rbrack = \frac{\ux (\uy \wedge \ux)}{t} \ \ \ \mathrm{and} \ \ \ \partial_{\uy} \lbrack \ux \wedge \uy \rbrack = (m-1) \ux,
\end{displaymath}
we obtain
\begin{align*}
\partial_{\uy} \lbrack K^+(\ux,\uy) \rbrack &= \ux \ \bigl( \partial_s \lbrack f(s,t) \rbrack + (m-1) \ g(s,t) + t \partial_t \lbrack g(s,t) \rbrack \bigr)\\
&+ \ux (\ux \wedge \uy) \ \left( \partial_s \lbrack g(s,t) \rbrack - \frac{1}{t} \partial_t \lbrack f(s,t) \rbrack \right),
\end{align*}
where we have used that $(\ux \wedge \uy)^2 = -t^2$. The right-hand side of the first equation of (\ref{system}) takes the form
\begin{displaymath}
a \ K^-(\ux,\uy) \ \ux = a \ f(-s,t) \ \ux - a \ g(-s,t) \ (\ux \wedge \uy) \ \ux.
\end{displaymath}
As 
\begin{displaymath}
\ux (\ux \wedge \uy) = \ux \ (\langle \ux,\uy \rangle + \ux \uy) = ( \langle \ux,\uy \rangle + \uy \ux) \ \ux = (\uy \wedge \ux) \ \ux = - (\ux \wedge \uy) \ \ux,
\end{displaymath}
we hence arrive at the following system for the functions $f$ and $g$:
\begin{align}\label{fg}
\begin{split}
\partial_s \lbrack f(s,t) \rbrack + t \partial_t \lbrack g(s,t)\rbrack + (m-1) \ g(s,t)  =  a \ f(-s,t)\\
\partial_s \lbrack g(s,t) \rbrack - \frac{1}{t} \partial_t \lbrack f(s,t) \rbrack  =  a \ g(-s,t).
\end{split}
\end{align}
The second equation of (\ref{system}) leads to the same system. 

We want to find solutions of the system (\ref{fg}) which are as close as possible to the kernel of the Clifford-Fourier transform given in formula (\ref{kernel C-F}). Therefore, we propose to find all solutions of the form
\begin{displaymath}
f(s,t)  =  \sum_{j=0}^k s^{k-j} \ f_j(t), \qquad g(s,t)  =  \sum_{j=0}^k s^{k-j} \ g_j(t)
\end{displaymath}
with $k \in \mathbb{N}$ a parameter. In other words, we want the solution to be polynomial in $s$, but do not prescribe the behavior of the $t$ variable.

Substituting this Ansatz in the system (\ref{fg}) yields 
\begin{align}
\label{fgj0}
(k-j+1) \ f_{j-1}(t) + (m-1) \ g_j(t) + t \ g_j'(t)  &=  a \ (-1)^{k-j} \ f_j(t)\\
\label{fgj1}
(k-j+1) \ g_{j-1}(t) -\frac{1}{t} \ f_j'(t)  &=  a \ (-1)^{k-j} \ g_j(t),
\end{align}
for $j = 0, \ldots, k$ and where $f_{-1}=g_{-1}=0$.

Let us first determine $f_0$ and $g_0$ from equations (\ref{fgj0}-\ref{fgj1}). Decoupling yields the following equation for $g_0$ :
\begin{displaymath}
t \ g_0''(t) + m \ g_0'(t) + t \ g_0(t) = 0
\end{displaymath}
from which we obtain (we want a solution which is not singular in $t=0$)
\begin{displaymath}
g_0(t) = c_0 \ \widetilde{J}_{(m-1)/2}(t), \qquad c_0 \in \mathbb{R}
\end{displaymath}
and thus also
\begin{displaymath}
f_0(t) = c_0 \ a \ (-1)^k \ \widetilde{J}_{(m-3)/2}(t).
\end{displaymath}
Subsequently, we determine $f_j$ and $g_j$ for $j=1,2,\ldots,k$ from equations (\ref{fgj0}-\ref{fgj1}).\\
We decouple the system by substituting (\ref{fgj1}) in the derivative of (\ref{fgj0}), yielding
\begin{equation}\label{eq9}
t g_j''(t) + m \ g_j'(t) + t \ g_j(t) = -(k-j+1) (k-j+2) \ t \ g_{j-2}(t) \ \ , \ \ j=2,3,\ldots,k.
\end{equation}
As we know $g_0$, the above differential equation yields the even $g's$ iteratively. Hereby we use the fact that $h(t)=\widetilde{J}_b(t)$ is a solution of the equation
\begin{displaymath}
t h''(t) + m \ h'(t) + t \ h(t) = (2b-m+1) \ t^{-b} \ J_{b+1}(t) 
\end{displaymath}
to determine a particular solution of (\ref{eq9}). In this way we obtain for $0 \leq \ell \leq \frac{k}{2}$ the following general solution for the even $g's$ (where at each step, we have excluded solutions singular in $t=0$):
\begin{displaymath}
g_{2 \ell}(t) = \sum_{i=0}^{\ell} c_{2i} \ \frac{1}{2^{\ell-i} (\ell-i)!} \ \frac{\Gamma(k+1-2i)}{\Gamma(k+1-2\ell)} \  \widetilde{J}_{(m-2\ell-1+2i)/2}(t) 
\end{displaymath}
with $c_{2i} \in \mathbb{R}$.

Next we determine the odd $g's$. Hereto we look for a differential equation for $g_1$. Substituting for $j=1$ the derivative of (\ref{fgj0}) in (\ref{fgj1}) leads to
\begin{align*}
t \ g_1''(t) + m \ g'_1(t) + t \ g_1(t) & = a \ (-1)^{k-1} \ k \ t \ g_0(t) - k \ f'_0(t)\\
\Leftrightarrow \  t \ g_1''(t) + m \ g'_1(t) + t \ g_1(t) & = 0,
\end{align*}
hence we find that $g_1(t) = c_1 \ \widetilde{J}_{(m-1)/2}(t)$ with $c_1 \in \mathbb{R}$. Now equation (\ref{eq9}) yields iteratively all odd $g's$ ($0 \leq \ell \leq \frac{k-1}{2}$):
\begin{displaymath}
g_{2 \ell+1}(t) = \sum_{i=0}^{\ell} c_{2i+1} \ \frac{1}{2^{\ell-i} (\ell-i)!} \ \frac{\Gamma(k-2i)}{\Gamma(k-2\ell)} \  \widetilde{J}_{(m-2\ell-1+2i)/2}(t) 
\end{displaymath}
with $c_{2i+1} \in \mathbb{R}$.

Next the $f_j$'s follow from equation (\ref{fgj0}):
\begin{align*}
f_{2 \ell}(t)  &= \sum_{i=0}^{\ell}  \frac{1}{2^{\ell-i} (\ell-i)!} \ \frac{\Gamma(k+1-2i)}{\Gamma(k+1-2\ell)} \ \biggl( -c_{2i-1} (k-(2i-1)) + c_{2i} \ a  (-1)^k \biggr)\\
& \times \widetilde{J}_{(m-2\ell-3+2i)/2}(t), \qquad 0 \leq \ell \leq \frac{k}{2}, \qquad c_{-1} = 0\\
f_{2 \ell+1}(t)  &= \sum_{i=0}^{\ell}  \frac{1}{2^{\ell-i} (\ell-i)!} \ \frac{\Gamma(k-2i)}{\Gamma(k-2\ell)} \ \left( -c_{2i} (k-2i) + c_{2i+1} \ a  (-1)^{k-1} \right)\\
& \times  \widetilde{J}_{(m-2\ell-3+2i)/2}(t), \qquad 0 \leq \ell \leq \frac{k-1}{2}. 
\end{align*}

If we change the summation order and renumber the integration coefficients as $c_{j} \rightarrow c_{k-j}$, we can summarize the class of solutions we have obtained as follows. We have
\[
K^+(\ux,\uy)  =  f(s,t) + (\ux \wedge \uy) \ g(s,t)
\]
with
\begin{align*}
f(s,t) &=  - \sum_{i = 1}^{k} c_{i} \  \sum_{\ell =0}^{\left\lfloor  \frac{i-1}{2} \right\rfloor} s^{i-2 \ell-1} \ \frac{1}{2^{\ell} \ell!} \ \frac{\Gamma(i+1)}{\Gamma(i-2\ell)} \ \widetilde{J}_{(m-2\ell-3)/2}(t)\\
& \ \ + a \ \sum_{i=0}^k c_i \ (-1)^i \ \sum_{\ell =0}^{\left\lfloor  \frac{i}{2}\right\rfloor} s^{i-2\ell} \ \frac{1}{2^{\ell} \ell!} \ \frac{\Gamma(i+1)}{\Gamma(i-2\ell+1)} \ \widetilde{J}_{(m-2\ell-3)/2}(t)\\
&= \sum_{i=0}^{k} c_i \ f_m^i(s,t)
\end{align*} 
where $f_m^i$ is independent of $k$. The function $g$ is given by
\begin{align*}
g(s,t) &=  \sum_{i=0}^k c_i \ \sum_{\ell=0}^{\left\lfloor  \frac{i}{2} \right\rfloor} s^{i-2\ell} \ \frac{1}{2^{\ell} \ell!} \ \frac{\Gamma(i+1)}{\Gamma(i+1-2\ell)} \ \widetilde{J}_{(m-2\ell-1)/2}(t)\\
&=  \sum_{i=0}^{k} c_i \ g_m^i(s,t)
\end{align*}
with again $g_m^i$ independent of $k$.

We want that the Bessel functions in this class of solutions are of order $\geq -1/2$. Therefore, we will restrict ourselves to the case where $k \leq m-2$.

%%%%%%%%%%%%%%%%%%%%%%%%%%%%%%%%%%%%%%%%%%%%%%%%%%%%%%%%%%%%%%%%%%%%%%%%%%%%%%%%%%%%%%%%%%%%%%%%%%%%%%%
\subsection{Recursion relations}

In this section we put
\begin{equation}\label{kernel even m}
K^{i}_{+ , m}(\ux,\uy)  =  \tilde{f}_m^i(s,t) +  \hat{f}_m^i(s,t) + (\ux \wedge \uy) \ g_m^i(s,t), \quad i=0,1,2, \ldots , m-2
\end{equation}
with
\begin{align*}
\tilde{f}_m^i(s,t) &=  -  \sqrt{\frac{\pi}{2}} \ \sum_{\ell =0}^{\left\lfloor  \frac{i-1}{2} \right\rfloor} s^{i-1-2 \ell} \ \frac{1}{2^{\ell} \ell!} \frac{\Gamma(i+1)}{\Gamma(i-2\ell)} \ \widetilde{J}_{(m-2\ell-3)/2}(t), \ \ i \geq 1\\
\hat{f}_m^i (s,t) &=  (-1)^{m/2+i} \   \sqrt{\frac{\pi}{2}} \ \sum_{\ell =0}^{\left\lfloor  \frac{i}{2} \right\rfloor} s^{i-2 \ell} \ \frac{1}{2^{\ell} \ell!} \frac{\Gamma(i+1)}{\Gamma(i+1-2\ell)} \ \widetilde{J}_{(m-2\ell-3)/2}(t) , \ \ i \geq 0\\
g_m^i(s,t) &=   \sqrt{\frac{\pi}{2}} \ \sum_{\ell =0}^{\left\lfloor  \frac{i}{2} \right\rfloor} s^{i-2 \ell} \ \frac{1}{2^{\ell} \ell!} \frac{\Gamma(i+1)}{\Gamma(i+1-2\ell)} \ \widetilde{J}_{(m-2\ell-1)/2}(t) , \ \  i \geq 0.
\end{align*}
Note that
\begin{itemize}
\item $K^{i}_{+ , m}$ ($i=0,1,2, \ldots , m-2$) is the solution with $c_i = \sqrt{\frac{\pi}{2}}$ and $c_0= \ldots = c_{i-1} = c_{i+1} = \ldots = c_{m-2} = 0$.
\item the Fourier-Bessel kernel (see (\ref{kernel F-B})) is obtained for $i=0$. Hence we put $K_{+,m}^{0}(\ux,\uy) = K_{+,m}^{\mathrm{Bessel}} (\ux,\uy)$. 
\item the even dimensional Clifford-Fourier kernel (see (\ref{kernel C-F})) is obtained, up to a minus sign, for $i = \frac{m}{2}-1$. Hence we denote $K_{+,m}^{m/2-1}(\ux,\uy) = K_{+,m}^{\mathrm{CF}} (\ux,\uy)$.
\item as $m$ is even, the solution is given in terms of Bessel functions of order $n+\frac{1}{2}$ with $n \in \mathbb{N}$.
\end{itemize}
We can arrange all the kernels as in the scheme below. Observe that at each step in the dimension, two new kernels appear ($K_{+,m}^0$ and $K_{+,m}^{m-2}$) corresponding to the Fourier-Bessel kernel and its inverse (as we will show in Theorem \ref{TFonS}). The other kernels at a given step in the dimension ($K_{+,m}^i$, $i=1,2,\ldots,m-3$) follow from the previous dimension $m-2$ by a suitable action of a differential operator as is explained in the following proposition. The middle line in the diagram corresponds with the Clifford-Fourier kernel.
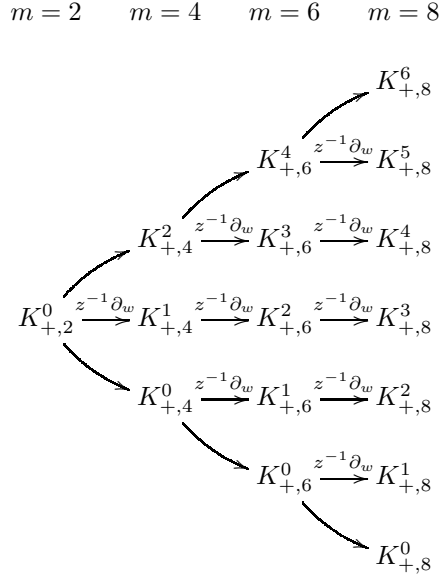
\[
\xymatrix@=12pt{m=2&m=4&m=6&m=8\\
&&&K_{+,8}^{6}\\
&&K_{+,6}^{4}\ar@/^/[ur] \ar[r]^-{z^{-1}\partial_{w}}&K_{+,8}^{5}\\
&K_{+,4}^{2} \ar@/^/[ur]\ar[r]^-{z^{-1}\partial_{w}}&K_{+,6}^{3}\ar[r]^-{z^{-1}\partial_{w}}&K_{+,8}^{4}\\
K_{+,2}^{0} \ar@/^/[ur] \ar@/_/[dr] \ar[r]^-{z^{-1}\partial_{w}}&K_{+,4}^{1} \ar[r]^-{z^{-1}\partial_{w}}&K_{+,6}^{2} \ar[r]^-{z^{-1}\partial_{w}}&K_{+,8}^{3}\\
&K_{+,4}^{0} \ar@/_/[dr] \ar[r]^-{z^{-1}\partial_{w}}&K_{+,6}^{1}\ar[r]^-{z^{-1}\partial_{w}}&K_{+,8}^{2}\\
&&K_{+,6}^{0} \ar@/_/[dr] \ar[r]^-{z^{-1}\partial_{w}}&K_{+,8}^{1}\\
&&&K_{+,8}^{0}
}
\]

\begin{proposition}~\\
A) For $1 \leq i \leq \frac{m}{2} - 1$ (lower half of the triangle in the above scheme) we have the following recursion relations :
\begin{align}
\label{rel1}
\tilde{f}_{m+2}^{i+1}(s,t) &= \frac{i+1}{i} \ z^{-1} \partial_w \tilde{f}_m^i(s,t)\\
\label{rel2}
\hat{f}_{m+2}^{i+1}(s,t) &=   z^{-1} \partial_w \hat{f}_m^i(s,t) \\
\label{rel3}
g^{i+1}_{m+2}(s,t) &=  - \frac{1}{i+1} \ z^{-1} \partial_w \tilde{f}_{m+2}^{i+1}(s,t)
\end{align}
with starting values given by the Fourier-Bessel kernel:
\begin{align}\label{rec1}
\begin{split}
\hat{f}^1_{m+2}(s,t) &=  z^{-1} \partial_w \hat{f}_m^0(s,t)\\
\tilde{f}^1_{m+2}(s,t)  &=  (-1)^{m/2-1} \ s^{-1} \ \hat{f}^1_{m+2}(s,t)\\
g^1_{m+2}(s,t)  &=  - z^{-1} \partial_w \tilde{f}^1_{m+2}(s,t)  .
\end{split}
\end{align}
B) For $\frac{m}{2} \leq i \leq m-2$ (upper half of the triangle in the above scheme) we have again the recursion relations (\ref{rel1}-\ref{rel3}), but now we start from the kernel $K_{+,m}^{m-2}$ :

\begin{align}\label{rec2}
\begin{split}
\tilde{f}_{m+2}^{m-1}(s,t)  &=  \frac{m-1}{m-2} \ z^{-1} \partial_w \tilde{f}_m^{m-2}(s,t)\\
\hat{f}_{m+2}^{m-1}(s,t)  &=   z^{-1} \partial_w \hat{f}_m^{m-2}(s,t) \\
g^{m-1}_{m+2}(s,t)  &=  - \frac{1}{m-1} \ z^{-1} \partial_w \tilde{f}_{m+2}^{m-1}(s,t)  .
\end{split}
\end{align}
Here, the notations $z=|\ux |  |\uy|$ and $w = \langle \underline{\xi} , \underline{\eta} \rangle$ ($\ux = |\ux | \ \underline{\xi}$, $\uy = |\uy| \ \underline{\eta}$) are used. Hence the transformation formulas $s = z w $ and $t=z \sqrt{1-w^2}$ hold.
\end{proposition}
\begin{proof}
The proof is carried out by induction on the dimension $m$ and is based on the properties:
\begin{displaymath}
z^{-1}  \partial_w \ \widetilde{J}_{\alpha}(t) = s \ \widetilde{J}_{\alpha + 1}(t) \qquad \mathrm{and} \qquad z^{-1}  \partial_w s^{\alpha} = \alpha \ s^{\alpha - 1} .
\end{displaymath}
By means of $\hat{f}_m^0(s,t) = (-1)^{m/2} \ \sqrt{\frac{\pi}{2}} \ \widetilde{J}_{(m-3)/2}(t)$ and making a distinction between $m=4p$ and $m=4p+2$, it is easy to check that the formulas (\ref{rec1}) are correct.\\
Next, let us check (\ref{rel1}) for the case i) $m=4p+4$ and $i=2j$. We have consecutively 
\begin{align*}
z^{-1} \partial_w \ \tilde{f}^{2j}_{4p+4}(s,t) &=  z^{-1} \partial_w \biggl( - \sqrt{\frac{\pi}{2}} \ \sum_{\ell=0}^{j-1} s^{2j-2\ell-1} \frac{1}{2^{\ell} \ell!} \frac{\Gamma (2j+1)}{\Gamma (2j -2 \ell)} \ \widetilde{J}_{2p+1/2-\ell}(t) \biggr)\\
&=  - \sqrt{\frac{\pi}{2}} \ \sum_{\ell=0}^{j-1} s^{2j-2\ell-2} \frac{1}{2^{\ell} \ell!} \frac{\Gamma (2j+1)}{\Gamma (2j -2 \ell-1)} \ \widetilde{J}_{2p+1/2-\ell}(t)\\
& \ \ - \sqrt{\frac{\pi}{2}} \ \sum_{\ell=0}^{j-1} s^{2j-2\ell} \frac{1}{2^{\ell} \ell!} \frac{\Gamma (2j+1)}{\Gamma (2j -2 \ell)} \ \widetilde{J}_{2p+3/2-\ell}(t)\\
&=  - \sqrt{\frac{\pi}{2}} \ \sum_{\ell=1}^{j} s^{2j-2\ell} \frac{2 \ell}{2^{\ell} \ell!} \frac{\Gamma (2j+1)}{\Gamma (2j -2 \ell+1)} \ \widetilde{J}_{2p+3/2-\ell}(t)\\
& \ \  - \sqrt{\frac{\pi}{2}} \ \sum_{\ell=0}^{j-1} s^{2j-2\ell} \frac{1}{2^{\ell} \ell!} \frac{\Gamma (2j+1)}{\Gamma (2j -2 \ell+1)} \ (2j-2\ell) \ \widetilde{J}_{2p+3/2-\ell}(t)
%\end{align*}
%\begin{align*}
\\
&=  - \sqrt{\frac{\pi}{2}} \ \sum_{\ell=1}^{j-1}  s^{2j-2\ell} \frac{1}{2^{\ell} \ell!} \frac{\Gamma (2j+1)}{\Gamma (2j -2 \ell+1)} \ (2j) \  \widetilde{J}_{2p+3/2- \ell}(t)\\
& \ \  - \sqrt{\frac{\pi}{2}} \ \frac{2j}{2^j j!} \ \Gamma (2j+1) \ \widetilde{J}_{2p-j+3/2}(t) - \sqrt{\frac{\pi}{2}} \ s^{2j} \ (2j) \ \widetilde{J}_{2p+3/2}(t)\\
&=  \frac{2j}{2j+1} \ \biggl( - \sqrt{\frac{\pi}{2}} \ \sum_{\ell=0}^{j} s^{2j-2\ell} \frac{1}{2^{\ell} \ell!} \frac{\Gamma (2j+2)}{\Gamma (2j -2 \ell+1)} \ \widetilde{J}_{2p+3/2-\ell}(t) \biggr)\\
&=  \frac{2j}{2j+1} \ \tilde{f}_{4p+6}^{2j+1}(s,t).
\end{align*}
The other cases: ii) $m=4p+4$, $i=2j+1$; iii) $m=4p+2$, $i=2j$; iv) $m=4p+2$, $i=2j+1$ are treated similarly.\\
The proof of formulas (\ref{rel2}) and (\ref{rel3}) runs along the same lines.\\
The formulas (\ref{rec2}) are an application of (\ref{rel1}-\ref{rel3}).
\end{proof}
%%%%%%%%%%%%%%%%%%%%%%%%%%%%%%%%%%%%%%%%%%%%%%%%%%%%%%%%%%%%%%%%%%%%%%%%%%%%%%%%%%%%%%%%%%%%%%
%%%%%%%%%%%%%%%%%%%%%%%%%%%%%%%%%%%%%%%%%%%%%%%%%%%%%%%%%%%%%%%%%%%%%%%%%%%%%%%%%%%%%%%%%%%%%%%%%%%%%
\subsection{Series expansion of $K_{+,m}^i$}
In this subsection we determine the series expansion in terms of Bessel functions and Gegenbauer polynomials of the kernels:
\[
K_{+,m}^{i}(\ux,\uy) = \tilde{f}_m^i + \hat{f}_m^i + (\ux \wedge \uy) \ g_m^i, \quad i = 0,1, 2, \ldots , m-2.
\]

\begin{theorem} 
\label{SeriesEven}
The following series expansions hold:\vspace{0,2cm}\\
\emph{\underline{Case 1}: $i$ even ($i=0,2,\ldots,m-4,m-2$)}
\begin{align*}
\tilde{f}_m^i(w,z) &=  - i \ \left( \frac{m}{2}-2 \right)! \ 2^{m/2-2}  \  \sum_{j=0}^{\infty} (4j+m) \ \frac{(2j+i-1)!!}{(2j+m-i-1)!!}\\
 & \times  z^{-m/2+1} \ J_{2j+m/2}(z) \ C_{2j+1}^{m/2-1}(w)\\
 \hat{f}_m^i(w,z) &=  (-1)^{m/2} \ \left( \frac{m}{2}-2 \right)! \ 2^{m/2-1} \ \sum_{j=0}^{\infty} \left( 2j + \frac{m}{2}-1 \right) \ \frac{(2j+i-1)!!}{(2j-i+m-3)!!}\\
& \times z^{1-m/2} \ J_{2j+m/2-1}(z) \ C_{2j}^{m/2-1}(w)\\
g_m^i(w,z) &=   \left( \frac{m}{2}-1 \right)! \ 2^{m/2-1} \ \sum_{j=0}^{\infty} (4j+m) \ \frac{(2j+i-1)!!}{(2j+m-i-1)!!}\\
& \times   z^{-m/2} \ J_{2j+m/2}(z) \ C_{2j}^{m/2}(w)
\end{align*}
%%%%%%%%%%%%%%%%%%%%%%%%%%%%%%%%%%%%%
\emph{\underline{Case 2}: $i$ odd ($i=1,3,\ldots, m-5, m-3$)}
\begin{align*}
\tilde{f}_m^i(w,z) &= - i \ \left( \frac{m}{2}-2 \right)! \ 2^{m/2-2}  \  \sum_{j=0}^{\infty} (4j+m-2) \ \frac{(2j+i-2)!!}{(2j+m-i-2)!!}\\
& \times   z^{-m/2+1} \ J_{2j+m/2-1}(z) \ C_{2j}^{m/2-1}(w)\\
\hat{f}_m^i(w,z) &= (-1)^{m/2+1} \  \left( \frac{m}{2}-2\right)! \ 2^{m/2-1} \ \sum_{j=0}^{\infty} \left( 2j + \frac{m}{2} \right) \ \frac{(2j+i)!!}{(2j+m-i-2)!!}\\
& \times  z^{1-m/2} \ J_{2j+m/2}(z) \ C_{2j+1}^{m/2-1}(w)\\
%\end{align*}
%\begin{align*}
g_m^i(w,z) &=   \left( \frac{m}{2}-1 \right)! \ 2^{m/2-1} \ \sum_{j=0}^{\infty}(4j+m+2) \  \frac{(2j+i)!!}{(2j+m-i)!!}\\
& \times    z^{-m/2} \ J_{2j+m/2+1}(z) \ C_{2j+1}^{m/2}(w) .
\end{align*}
\end{theorem}

\begin{proof}
We prove this theorem by induction on the dimension $m$.\\
For $m=2$, the kernel $K_{+,2}^0$ is the Clifford-Fourier transform for which the series expansion is derived in \cite{DBXu}.\\
Suppose we know the series expansion for all kernels in dimension $m-2$ and lower, then we prove that we can obtain all series expansions in dimension $m$.\\
The terms $\tilde{f}_m^i$, $\hat{f}_m^i$, $g_m^i$, $i=1,2,\ldots,m-3$ are easy; they follow immediately by action of $z^{-1} \partial_w$ on the terms in dimension $m-2$ as is indicated for $g_m^i$ in the diagram below (where we have omitted the factor $\sqrt{\frac{\pi}{2}}$). 

\[
\xymatrix@=12pt{m=2&m=4&m=6&m=8\\
&&& \ldots\\
&&s^{4}\widetilde{J}_{5/2}+ 6 s^{2} \widetilde{J}_{3/2}+ 3 \widetilde{J}_{1/2}\ar[r]^-{z^{-1}\partial_{w}}& \ldots\\
&s^{2}\widetilde{J}_{3/2} + \widetilde{J}_{1/2}\ar[r]^-{z^{-1}\partial_{w}}&s^{3 }\widetilde{J}_{5/2}+ 3 s \widetilde{J}_{3/2}\ar[r]^-{z^{-1}\partial_{w}}& \ldots\\
\widetilde{J}_{1/2} \ar[r]^-{z^{-1}\partial_{w}}&s \widetilde{J}_{3/2} \ar@/^/[d]^{\frac{1}{s}} \ar[r]^-{z^{-1}\partial_{w}}&s^{2}\widetilde{J}_{5/2}+ \widetilde{J}_{3/2}\ar[r]^-{z^{-1}\partial_{w}}& \ldots\\
&\widetilde{J}_{3/2}\ar[r]^-{z^{-1}\partial_{w}}&s\widetilde{J}_{5/2} \ar@/^/[d]^{\frac{1}{s}}\ar[r]^-{z^{-1}\partial_{w}}& \ldots \\
&&\widetilde{J}_{5/2}\ar[r]^-{z^{-1}\partial_{w}}&s \widetilde{J}_{7/2}\ar@/^/[d]^{\frac{1}{s}}  \\
&&&\widetilde{J}_{7/2}
}
\]

For example, let us consider the term $\tilde{f}_m^i$ with $i$ even ($i=2,4,\ldots,m-4$). By means of (\ref{rel1}) and the already known expansion of $\tilde{f}_{m-2}^{i-1}$ we have that

\begin{align*}
\tilde{f}_m^i & = \frac{i}{i-1} \ z^{-1}\partial_w \tilde{f}_{m-2}^{i-1}\\
&= \frac{i}{i-1} \ z^{-1}\partial_w  \biggl( - (i-1) \ \left( \frac{m}{2}-3 \right)! \ 2^{m/2-3}  \  \sum_{j=0}^{\infty} (4j+m-4)\\
& \ \ \ \ \ \times \frac{(2j+i-3)!!}{(2j+m-i-3)!!} \ z^{-m/2+2} \ J_{2j+m/2-2}(z) \ C_{2j}^{m/2-2}(w) \biggr)\\
&=  - i \ \left( \frac{m}{2}-2 \right)! \ 2^{m/2-2}  \  \sum_{j=1}^{\infty} (4j+m-4) \  \frac{(2j+i-3)!!}{(2j+m-i-3)!!}\\
& \ \ \ \ \ \times z^{-m/2+1} \ J_{2j+m/2-2}(z) \ C_{2j-1}^{m/2-1}(w)\\
%\end{align*}
%\begin{align*}
&= - i \ \left( \frac{m}{2}-2 \right)! \ 2^{m/2-2}  \  \sum_{j=0}^{\infty} (4j+m) \  \frac{(2j+i-1)!!}{(2j+m-i-1)!!} \ z^{-m/2+1}\\
& \ \ \ \ \ \times J_{2j+m/2}(z) \ C_{2j+1}^{m/2-1}(w),
\end{align*}
where in the third step we have used $\frac{d}{dw} \left( C_k^{\lambda}(w) \right) = 2 \lambda \ C_{k-1}^{\lambda + 1}(w)$.

So we only need to find the series expansions of $\hat{f}_m^0$, $g_m^0$, $\tilde{f}_m^{m-2}$, $\hat{f}_m^{m-2}$ and $g_m^{m-2}$. Let us first treat the terms $g_m^0$ and $g_m^{m-2}$ and afterwards describe the similar procedure for $\hat{f}_m^0$, $\tilde{f}_m^{m-2}$ and $\hat{f}_m ^{m-2}$.

The kernel $g_m^0$ can be derived from $g_m^1$ because
\begin{displaymath}
g_m^0 = \frac{1}{s} \ g_m^1.
\end{displaymath}
We have by using the series expansion of $g_m^1$ and the Bessel identity (\ref{Besselid}) that 
\begin{align*}
g_m^0 &= \frac{1}{wz} \ \left( \frac{m}{2}-1 \right)! \ 2^{m/2-1} \ \sum_{j=0}^{\infty} (4j+m+2) \ \frac{(2j+1)!!}{(2j+m-1)!!} \ z^{-m/2} \\
& \ \  J_{2j+m/2+1}(z) \  C_{2j+1}^{m/2}(w)\\
&=  \frac{\left( \frac{m}{2}-1 \right)! \ 2^{m/2-1}}{w} \biggl( \sum_{j=0}^{\infty} \frac{(2j+1)!!}{(2j+m-1)!!} \ z^{-m/2} \ J_{2j+m/2}(z) \ C_{2j+1}^{m/2}(w)\\
& \ \  + \sum_{j=0}^{\infty} \frac{(2j+1)!!}{(2j+m-1)!!} \ z^{-m/2} \ J_{2j+m/2+2}(z) \ C_{2j+1}^{m/2}(w) \biggr) .
\end{align*}
Next, executing the substitution $j=j'-1$ in the second term and applying the formula (\ref{Geg2}) yields the desired result:
\begin{align*}
g_m^0 &=  \frac{\left( \frac{m}{2}-1 \right)! \ 2^{m/2-1}}{w} \biggl( \sum_{j=0}^{\infty} \frac{(2j+1)!!}{(2j+m-1)!!} \ z^{-m/2} \ J_{2j+m/2}(z) \ C_{2j+1}^{m/2}(w)\\
& \ \  + \sum_{j=1}^{\infty} \frac{(2j-1)!!}{(2j+m-3)!!} \ z^{-m/2} \ J_{2j+m/2}(z) \ C_{2j-1}^{m/2}(w) \biggr)\\
&=  \frac{\left( \frac{m}{2}-1 \right)! \ 2^{m/2-1}}{w} \ \sum_{j=0}^{\infty} \frac{(2j-1)!!}{(2j+m-1)!!} \ z^{-m/2} \ J_{2j+m/2}(z)\\
& \ \  \biggl( (2j+1) \  C_{2j+1}^{m/2}(w) +  (2j+m-1) \  C_{2j-1}^{m/2}(w) \biggr)\\
&=  \left( \frac{m}{2}-1 \right)! \ 2^{m/2-1} \ \sum_{j=0}^{\infty} (4j+m) \ \frac{(2j-1)!!}{(2j+m-1)!!} \ z^{-m/2} \ J_{2j+m/2}(z) \ C_{2j}^{m/2}(w) .
\end{align*}

The series expansion for the term $g_{m}^{m-2}$ follows from the following relation:
\begin{displaymath}
g_m^{m-2} = s \ g_m^{m-3} + (m-3) \ g_{m-2}^{m-4} 
\end{displaymath}
which can be observed from the diagram with the explicit expressions for $g_m^i$ ($m=2,4,6,8$) and proved by a direct calculation.

We now use the already known expansions of $g_m^{m-3}$ and $g_{m-2}^{m-4}$, combined with use of the formula (\ref{Geg2}) on the first term and (\ref{Geg1}) on the second term, yielding
\begin{align*}
&  s \ g_m^{m-3} + (m-3) \ g_{m-2}^{m-4} \\
&=   \left( \frac{m}{2}-1 \right)! \ 2^{m/2-1} \ \sum_{j=0}^{\infty} (4j+m+2) \ \frac{(2j+m-3)!!}{(2j+3)!!} \ z^{-m/2+1} \ J_{2j+m/2+1}(z) \ w\\
& \ \    C_{2j+1}^{m/2}(w) + (m-3) \ \left( \frac{m}{2}-2 \right)! \ 2^{m/2-2} \ \sum_{j=0}^{\infty} (4j+m-2) \ \frac{(2j+m-5)!!}{(2j+1)!!} \ z^{-m/2+1} \\
& \ \  J_{2j+m/2-1}(z) \  C_{2j}^{m/2-1}(w)\\
&=  \left( \frac{m}{2}-1 \right)! \ 2^{m/2-1} \  \biggl( \sum_{j=0}^{\infty} \frac{(2j+m-3)!!}{(2j+3)!!} \ z^{-m/2+1} \ J_{2j+m/2+1}(z) \ (2j+2) \  C_{2j+2}^{m/2}(w)\\
& \ \  + \sum_{j=0}^{\infty} \frac{(2j+m-3)!!}{(2j+3)!!} \ z^{-m/2+1} \ J_{2j+m/2+1}(z) \ (2j+m) \  C_{2j}^{m/2}(w)\\
& \ \  + (m-3) \ \sum_{j=0}^{\infty} \frac{(2j+m-5)!!}{(2j+1)!!} \ z^{-m/2+1} \ J_{2j+m/2-1}(z)  \  C_{2j}^{m/2}(w)\\
& \ \ - (m-3) \ \sum_{j=0}^{\infty} \frac{(2j+m-5)!!}{(2j+1)!!} \ z^{-m/2+1} \ J_{2j+m/2-1}(z)  \  C_{2j-2}^{m/2}(w) \biggr) .
\end{align*}
Consecutively, we execute in the first term the substitution $j=j'-1$ and in the last term the substitution $j=j'+1$, after which we can collect the first term and the third one and similarly the second term and the last one. In this way we arrive at:
\begin{align*}
&  s \ g_m^{m-3} + (m-3) \ g_{m-2}^{m-4} \\
&=   \left( \frac{m}{2}-1 \right)! \ 2^{m/2-1} \ \biggl( \sum_{j=0}^{\infty} \frac{(2j+m-3)!!}{(2j+1)!!} \ z^{-m/2+1} \ J_{2j+m/2-1}(z) \  C_{2j}^{m/2}(w)\\
& \ \ + \sum_{j=0}^{\infty} \frac{(2j+m-3)!!}{(2j+1)!!} \ z^{-m/2+1} \ J_{2j+m/2+1}(z) \  C_{2j}^{m/2}(w) \biggr) .
\end{align*}
Moreover, using the Bessel identity (\ref{Besselid}) yields the desired series expansion for $g_m^{m-2}$ :
\begin{eqnarray*}
& & s \ g_m^{m-3} + (m-3) \ g_{m-2}^{m-4} \\
& = &  \left( \frac{m}{2}-1 \right)! \ 2^{m/2-1} \ \sum_{j=0}^{\infty} (4j+m) \ \frac{(2j+m-3)!!}{(2j+1)!!} \ z^{-m/2} \ J_{2j+m/2}(z) \  C_{2j}^{m/2}(w)  .
\end{eqnarray*}

The series expansions of $\hat{f}_m^0$, $\hat{f}_m^{m-2}$ and $\tilde{f}_m^{m-2}$ are derived in an analogous manner.
In case of $\hat{f}_m^0$ the expansion follows from the observation that
\begin{displaymath}
\hat{f}_m^0 = - (-1)^{m/2} \tilde{f}_m^1,
\end{displaymath}
while for $\hat{f}_m^{m-2}$ we must use the formula
\begin{displaymath}
\hat{f}_m^{m-2} = -s \ \hat{f}_m^{m-3} - (m-3) \ \hat{f}_{m-2}^{m-4}  .
\end{displaymath}
Finally, the series expansion for $\tilde{f}_m^{m-2}$ follows from the one of $\tilde{f}_m^{m-3}$ and $\tilde{f}_{m-2}^{m-4}$ by applying in a similar way as above the relation:
\begin{displaymath}
\tilde{f}_m^{m-2} = \frac{m-2}{m-3} \ s \ \tilde{f}_m^{m-3} + (m-2) \ \tilde{f}_{m-2}^{m-4}, 
\end{displaymath} 
which can again be proven by direct computation.
\end{proof}
%%%%%%%%%%%%%%%%%%%%%%%%%%%%%%%%%%%%%%%%%%%%%%%%%%%%%%%%%%%%%%%%%%%%%%%%%%%%%%%%%%%%%%%%%%%%%%%%%%%%%%%%%%%%%%%%%%%%%%%%%
\subsection{Eigenvalues of new class of Clifford-Fourier transforms}
\label{eigvalseven}

We will now calculate the action of the new Clifford-Fourier transforms
\begin{displaymath}
\mathcal{F}^{i}_{+,m} \lbrack f(\ux) \rbrack (\uy) = \frac{1}{(2\pi)^{m/2}} \ \int_{\mathbb{R}^m} K^{i}_{+,m}(\ux,\uy) \ f(\ux) \ dV(\ux), \qquad i=0,1,2,\ldots,m-2
\end{displaymath}
on the basis $\lbrace \psi_{j,k,\ell} \rbrace$. 

From the previous subsection we observe that the new Clifford-Fourier kernels $K_{+,m}^{i}$ are of the structure (\ref{structure kernel}). The action of the Clifford-Fourier transforms $\mathcal{F}^{i}_{+,m}$ on the basis $\lbrace \psi_{j,k,\ell} \rbrace$ can hence be determined by substituting the corresponding coefficients $\alpha_k$ and $\beta_k$ in the equations (see (\ref{eigenvalue eq})):
\begin{align}\label{eigenvalue eq scaling}
\begin{split}
\mathcal{F}^{i}_{+,m} \lbrack \psi_{2p,k,\ell} \rbrack (\uy) &=  \frac{2^{1-m/2}}{\Gamma \left( \frac{m}{2} \right)} \ \left( \frac{\frac{m}{2}-1}{\frac{m}{2}-1+k} \ \alpha_k - \frac{k}{2(\frac{m}{2}-1+k)} \ \beta_k \right) \ (-1)^p \ \psi_{2p,k,\ell}(\uy)\\
\mathcal{F}^{i}_{+,m} \lbrack \psi_{2p+1,k,\ell} \rbrack (\uy) &=  \frac{2^{1-m/2}}{\Gamma \left( \frac{m}{2} \right)} \ \left( \frac{\frac{m}{2}-1}{\frac{m}{2}+k} \ \alpha_{k+1} + \frac{k+m-1}{2(\frac{m}{2}+k)} \ \beta_{k+1} \right) \ (-1)^p \ \psi_{2p+1,k,\ell}(\uy) . 
\end{split}
\end{align}

This yields the following result:
\begin{theorem}
\label{EigValsEven}
In case of $m$ even, the Clifford-Fourier transforms $\mathcal{F}^{i}_{+,m}$ act as follows on the basis $\lbrace \psi_{j,k,\ell} \rbrace$ of ${\mathcal S}(\mathbb{R}^m) \otimes \cC l_{0,m}$:\vspace{0,2cm}\\
\emph{\underline{Case 1}: $i$ even} ($i=0,2,\ldots,m-2$).\vspace{0,1cm}\\
\underline{a) $k$: even}
\begin{align*}
\mathcal{F}^{i}_{+,m} \lbrack \psi_{2p,k,\ell} \rbrack (\uy) &= (-1)^{m/2} \ \frac{(k+i-1)!!}{(k-i+m-3)!!} \ (-1)^p \ \psi_{2p,k,\ell}(\uy)\\
\mathcal{F}^{i}_{+,m} \lbrack \psi_{2p+1,k,\ell} \rbrack (\uy) &=   \frac{(k+i-1)!!}{(k+m-i-3)!!} \ (-1)^p \ \psi_{2p+1,k,\ell}(\uy)
\end{align*}
\underline{b) $k$: odd}
\begin{align*}
\mathcal{F}^{i}_{+,m} \lbrack \psi_{2p,k,\ell} \rbrack (\uy) &=    \frac{(k+i)!!}{(k+m-i-2)!!} \ (-1)^{p+1} \ \psi_{2p,k,\ell}(\uy)\\
\mathcal{F}^{i}_{+,m} \lbrack \psi_{2p+1,k,\ell} \rbrack (\uy) &=  (-1)^{m/2} \  \frac{(k+i)!!}{(k-i+m-2)!!} \ (-1)^p \ \psi_{2p+1,k,\ell}(\uy)
\end{align*}
\emph{\underline{Case 2}: $i$ odd} ($i=1,3,\ldots,m-5,m-3$).\vspace{0,1cm}\\
\underline{a) $k$: even}
\begin{align*}
\mathcal{F}^{i}_{+,m} \lbrack \psi_{2p,k,\ell} \rbrack (\uy) &=   \frac{(k+i)!!}{(k+m-i-2)!!} \ (-1)^{p+1} \ \psi_{2p,k,\ell}(\uy)\\
\mathcal{F}^{i}_{+,m} \lbrack \psi_{2p+1,k,\ell} \rbrack (\uy) &=  (-1)^{m/2} \  \frac{(k+i)!!}{(k+m-i-2)!!} \ (-1)^{p+1} \ \psi_{2p+1,k,\ell}(\uy)
\end{align*}
\underline{b) $k$: odd}
\begin{align*}
\mathcal{F}^{i}_{+,m} \lbrack \psi_{2p,k,\ell} \rbrack (\uy) &=   (-1)^{m/2} \ \frac{(k+i-1)!!}{(k+m-i-3)!!} \ (-1)^{p+1} \ \psi_{2p,k,\ell}(\uy)\\
\mathcal{F}^{i}_{+,m} \lbrack \psi_{2p+1,k,\ell} \rbrack (\uy) &=   \frac{(k+i-1)!!}{(k+m-i-3)!!} \ (-1)^p \ \psi_{2p+1,k,\ell}(\uy).
\end{align*}
\end{theorem}

\begin{remark}
Putting $i=0$, we indeed obtain the eigenvalue equations of the Fourier-Bessel transform (see (\ref{eigenvalue F-B1}) and (\ref{eigenvalue F-B2})), while for $i=\frac{m}{2}-1$ the ones of the Clifford-Fourier transform appear (see (\ref{eigenvalue C-F})).

The result for the Fourier-Bessel transform was obtained in \cite{FourierBessel} using complicated integral identities for special functions. The present method is clearly more general and insightful.
\end{remark}
%%%%%%%%%%%%%%%%%%%%%%%%%%%%%%%%%%%%%%%%%%%%%%%%%%%%%%%%%%%%%%%%%%%%%%%%%%%%%%%%%%%%%%%%%%%%%%%%%%%%%%%ù
%%%%%%%%%%%%%%%%%%%%%%%%%%%%%%%%%%%%%%%%%%%%%%%%%%%%%%%%%%%%%%%%%%%%%%%%%%%%%%%%%%%%%%%%%%%%%%%%%%%%%%%%%%
\section{New Clifford-Fourier transforms: the case $m$ odd}
\setcounter{equation}{0}
\label{Odd}

\subsection{Parabivector-valued solutions of the Clifford-Fourier system}
We now want to solve the Clifford-Fourier system (\ref{C-F system}) in odd dimension:
\begin{align}\label{system2}
\begin{split}
\partial_{\uy} \lbrack K^+(\ux,\uy) \rbrack =  I a \ K^-(\ux,\uy) \ \ux\\
\lbrack K^+(\ux,\uy) \rbrack \partial_{\ux}  =  I a \ \uy \ K^-(\ux,\uy)
\end{split}
\end{align}
with $ a=(-1)^{(m+1)/2}$ and where $K^-(\ux,\uy) = \left( K^+(\ux,-\uy) \right)^c$, see (\ref{rel C-F kernel}).

We look for parabivector-valued solutions of the form
\begin{align*}
K^+(\ux,\uy)  &=  U(s,t)+ I \  V(s,t) + (\ux \wedge \uy) \ \lbrack Z(s,t) + I \ T(s,t) \rbrack\\
K^-(\ux,\uy)  &=  U(-s,t)- I \ V(-s,t) - (\ux \wedge \uy) \ \lbrack Z(-s,t) - I \ T(-s,t) \rbrack
\end{align*}
with again $s= \la \ux,\uy \ra$, $t=|\ux \wedge \uy|$ and $U$, $V$, $Z$ and $T$ real-valued functions, thus mimicking the form of the Clifford-Fourier kernel in formula (\ref{kernel C-F}).

Rewriting the system (\ref{system2}) in terms of $U$, $V$, $Z$ and $T$ yields
\begin{align}\label{systemUVZT}
\begin{split}
\partial_s \lbrack U(s,t) \rbrack + t \partial_t \lbrack Z(s,t) \rbrack + (m-1) \ Z(s,t)  =  a \ V(-s,t)\\
\partial_s \lbrack V(s,t) \rbrack + t \partial_t \lbrack T(s,t) \rbrack + (m-1) \ T(s,t)  =  a \ U(-s,t)\\
\partial_s \lbrack Z(s,t) \rbrack - \frac{1}{t} \partial_t \lbrack U(s,t) \rbrack    =  a \ T(-s,t)\\
\partial_s \lbrack T(s,t) \rbrack - \frac{1}{t} \partial_t \lbrack V(s,t) \rbrack    =  a \ Z(-s,t).
\end{split}
\end{align}
We are interested in solutions of the system (\ref{systemUVZT}) of the following form
\begin{align*}
U(s,t) &=  \sum_{j=0}^k s^{k-j} \ U_j(t) \quad , \quad V(s,t)  =  \sum_{j=0}^k s^{k-j} \ V_j(t)\\
Z(s,t) &=  \sum_{j=0}^k s^{k-j} \ Z_j(t) \quad , \quad T(s,t)  =  \sum_{j=0}^k s^{k-j} \ T_j(t)
\end{align*}
with $k \in \mathbb{N}$ a parameter. Using the same techniques as in the case $m$ even, one can explicitly determine such solutions. After lengthy computations, one finally arrives at the following general solution of the system (\ref{system2}):
\begin{displaymath}
K^+(\ux,\uy)  =  U(s,t) + I \ V(s,t) +  (\ux \wedge \uy) \ \lbrack Z(s,t) + I \ T(s,t) \rbrack
\end{displaymath}
with

\begin{align*}
U(s,t) + I \ V(s,t)  &=  - \sum_{i = 1}^{k} e_{i} \  \sum_{\ell =0}^{\left\lfloor  \frac{i-1}{2} \right\rfloor} s^{i-2 \ell-1} \ \frac{1}{2^{\ell} \ell!} \ \frac{\Gamma(i+1)}{\Gamma(i-2\ell)} \ \widetilde{J}_{(m-2\ell-3)/2}(t)\\
 & + a \  I \ \sum_{i = 0}^{k} e_{i}^c \ (-1)^{i}  \sum_{\ell =0}^{\left\lfloor  \frac{i}{2}\right\rfloor} s^{i-2\ell} \ \frac{1}{2^{\ell} \ell!} \ \frac{\Gamma(i+1)}{\Gamma(i-2\ell+1)} \ \widetilde{J}_{(m-2\ell-3)/2}(t)\\
Z(s,t) + I \ T(s,t)  &=  \sum_{i=0}^k e_i \ \sum_{\ell=0}^{\left\lfloor  \frac{i}{2} \right\rfloor} s^{i-2\ell} \ \frac{1}{2^{\ell} \ell!} \ \frac{\Gamma(i+1)}{\Gamma(i+1-2\ell)} \ \widetilde{J}_{(m-2\ell-1)/2}(t).
\end{align*} 
Here, $e_i \in \mathbb{C}$.

\begin{remark}
As $m$ is odd, the solution is given in terms of Bessel functions of integer order $n \in \mathbb{N}$. Again we take $k \leq m-2$, to ensure that the Bessel functions are of order $\geq 0$.
\end{remark}
%%%%%%%%%%%%%%%%%%%%%%%%%%%%%%%%%%%%%%%%%%%%%%%%%%%%%%%%%%%%%%%%%%%%%%%%%%%%%%%%%%%%%%%%%%%%%%%%%%%%%%%%%%%%%%
\subsection{Recursion relations}
In this section we put
\begin{equation}\label{kernel odd m}
K^{i}_{+ , m}(\ux,\uy)  =  e_i \ \tilde{f}_m^i(s,t) + I \ e_i^c \ \hat{f}_m^i(s,t) + (\ux \wedge \uy) \ e_i \ g_m^i(s,t),
\end{equation}
$i=0,1,2, \ldots , m-2$, with
\begin{align*}
\tilde{f}_m^i(s,t) &= -  \sum_{\ell =0}^{\left\lfloor  \frac{i-1}{2} \right\rfloor} s^{i-1-2 \ell} \ \frac{1}{2^{\ell} \ell!} \frac{\Gamma(i+1)}{\Gamma(i-2\ell)} \ \widetilde{J}_{(m-2\ell-3)/2}(t) , \ \ i \geq 1\\
\hat{f}_m^i (s,t) &=  (-1)^{(m+1)/2+i} \   \sum_{\ell =0}^{\left\lfloor  \frac{i}{2} \right\rfloor} s^{i-2 \ell} \ \frac{1}{2^{\ell} \ell!} \frac{\Gamma(i+1)}{\Gamma(i+1-2\ell)} \ \widetilde{J}_{(m-2\ell-3)/2}(t), \ \ i \geq 0\\
g_m^i(s,t) &=    \sum_{\ell =0}^{\left\lfloor  \frac{i}{2} \right\rfloor} s^{i-2 \ell} \ \frac{1}{2^{\ell} \ell!} \frac{\Gamma(i+1)}{\Gamma(i+1-2\ell)} \ \widetilde{J}_{(m-2\ell-1)/2}(t), \ \  i \geq 0 
\end{align*}
and $e_i \in \mathbb{C}$.\\
Note that $K_{+,m}^i$ is the $i$-th term in the general solution $K^+(\ux,\uy)$ of the previous subsection.

Similar to the even dimensional case, we can arrange all the kernels in the scheme below.
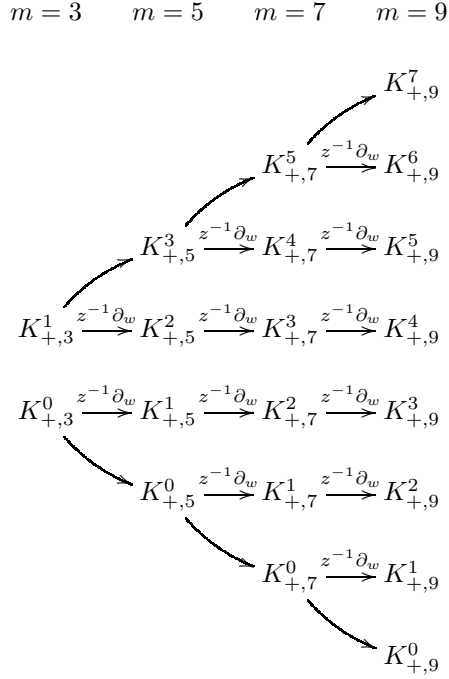
\[
\xymatrix@=13pt{m=3&m=5&m=7&m=9\\
&&&K_{+,9}^{7}\\
&&K_{+,7}^{5}\ar@/^/[ur] \ar[r]^-{z^{-1}\partial_{w}}&K_{+,9}^{6}\\
&K_{+,5}^{3} \ar@/^/[ur]\ar[r]^-{z^{-1}\partial_{w}}&K_{+,7}^{4}\ar[r]^-{z^{-1}\partial_{w}}&K_{+,9}^{5}\\
K_{+,3}^{1} \ar@/^/[ur]  \ar[r]^-{z^{-1}\partial_{w}}&K_{+,5}^{2} \ar[r]^-{z^{-1}\partial_{w}}&K_{+,7}^{3} \ar[r]^-{z^{-1}\partial_{w}}&K_{+,9}^{4}\\
K_{+,3}^{0}  \ar@/_/[dr] \ar[r]^-{z^{-1}\partial_{w}}&K_{+,5}^{1} \ar[r]^-{z^{-1}\partial_{w}}&K_{+,7}^{2} \ar[r]^-{z^{-1}\partial_{w}}&K_{+,9}^{3}\\
&K_{+,5}^{0} \ar@/_/[dr] \ar[r]^-{z^{-1}\partial_{w}}&K_{+,7}^{1}\ar[r]^-{z^{-1}\partial_{w}}&K_{+,9}^{2}\\
&&K_{+,7}^{0} \ar@/_/[dr] \ar[r]^-{z^{-1}\partial_{w}}&K_{+,9}^{1}\\
&&&K_{+,9}^{0}
}
\]

Observe that at each step in the dimension, two new kernels appear, namely $K_{+,m}^0$ and $K_{+,m}^{m-2}$. The other kernels at a given step in the dimension ($K_{+,m}^i$, $i=1,2,\ldots,m-3$) follow from the previous dimension $m-2$ by a suitable action of a differential operator as is explained in the following proposition. 

\begin{proposition}~\\
A) For $1 \leq i \leq \frac{m-1}{2} - 1$ (lower half of the triangle in the above scheme) we have the following recursion relations :
\begin{align}\label{recidem}
\begin{split}
\tilde{f}_{m+2}^{i+1}(s,t) &=  \frac{i+1}{i} \ z^{-1} \partial_w \tilde{f}_m^i(s,t)\\
\hat{f}_{m+2}^{i+1}(s,t) &=   z^{-1} \partial_w \hat{f}_m^i(s,t) \\
g^{i+1}_{m+2}(s,t) &=  - \frac{1}{i+1} \ z^{-1} \partial_w \tilde{f}_{m+2}^{i+1}(s,t)
\end{split}
\end{align}
with starting values given by the kernel $K_{+,m}^0$:
\begin{align*}
\hat{f}^1_{m+2}(s,t) &=  z^{-1} \partial_w \hat{f}_m^0(s,t)\\
\tilde{f}^1_{m+2}(s,t) &=  (-1)^{(m-1)/2} \ s^{-1} \ \hat{f}^1_{m+2}(s,t)\\
g^1_{m+2}(s,t) &=  - z^{-1} \partial_w \tilde{f}^1_{m+2}(s,t) .
\end{align*}
B) For $\frac{m-1}{2} \leq i \leq m-2$ (upper half of the triangle in the above scheme) we have again the recursion relations (\ref{recidem}), but now we start from the kernel $K_{+,m}^{m-2}$ :
\begin{align*}
\tilde{f}_{m+2}^{m-1}(s,t) &=  \frac{m-1}{m-2} \ z^{-1} \partial_w \tilde{f}_m^{m-2}(s,t)\\
\end{align*}
\begin{align*}
\hat{f}_{m+2}^{m-1}(s,t) &=   z^{-1} \partial_w \hat{f}_m^{m-2}(s,t) \\
g^{m-1}_{m+2}(s,t) &=  - \frac{1}{m-1} \ z^{-1} \partial_w \tilde{f}_{m+2}^{m-1}(s,t).
\end{align*}
\end{proposition}
\begin{proof}
The proof is similar to the even dimensional case.
\end{proof}
%%%%%%%%%%%%%%%%%%%%%%%%%%%%%%%%%%%%%%%%%%%%%%%%%%%%%%%%%%%%%%%%%%%%%%%%%%%%%%%%%%%%%%%%%%%%%%
%%%%%%%%%%%%%%%%%%%%%%%%%%%%%%%%%%%%%%%%%%%%%%%%%%%%%%%%%%%%%%%%%%%%%%%%%%%%%%%%%%%%%%%%%%%%%%%%%%%%%
\subsection{Series expansion of $K_{+,m}^i$}

Similar to the even dimensional case, we can expand the solutions in series in terms of Bessel functions and Gegenbauer polynomials.

\begin{theorem} 
\label{SeriesOdd}
The following series expansions hold:\vspace{0,2cm}\\
\emph{\underline{Case 1}: $i$ even ($i=0,2,\ldots,m-5,m-3$)}
\begin{align*}
\tilde{f}_m^i(w,z) &= - i \ 2^{m/2-2} \ \Gamma \left( \frac{m}{2}-1 \right)   \  \sum_{j=0}^{\infty} (4j+m) \ \frac{(2j+i-1)!!}{(2j+m-i-1)!!}\\
 & \times   z^{-m/2+1} \ J_{2j+m/2}(z) \ C_{2j+1}^{m/2-1}(w)\\
 \hat{f}_m^i(w,z) &=  (-1)^{(m+1)/2} \ 2^{m/2-1} \ \Gamma \left( \frac{m}{2}-1 \right) \ \sum_{j=0}^{\infty} \left( 2j + \frac{m}{2}-1 \right) \ \frac{(2j+i-1)!!}{(2j-i+m-3)!!}\\
  & \times  z^{1-m/2} \ J_{2j+m/2-1}(z) \ C_{2j}^{m/2-1}(w)\\
g_m^i(w,z) &=   2^{m/2-1} \ \Gamma \left( \frac{m}{2}\right) \ \sum_{j=0}^{\infty} (4j+m) \ \frac{(2j+i-1)!!}{(2j+m-i-1)!!}\\
& \times    z^{-m/2} \ J_{2j+m/2}(z) \ C_{2j}^{m/2}(w)
\end{align*}
%%%%%%%%%%%%%%%%%%%%%%%%%%%%%%%%%%%%%
\emph{\underline{Case 2}: $i$ odd ($i=1,3,\ldots, m-4, m-2$)}
\begin{align*}
\tilde{f}_m^i(w,z) &= - i \ 2^{m/2-2} \ \Gamma \left( \frac{m}{2}-1 \right)  \  \sum_{j=0}^{\infty} (4j+m-2) \ \frac{(2j+i-2)!!}{(2j+m-i-2)!!}\\
& \times   z^{-m/2+1} \ J_{2j+m/2-1}(z) \ C_{2j}^{m/2-1}(w)\\
\hat{f}_m^i(w,z) &=  - (-1)^{(m+1)/2} \ 2^{m/2-1} \ \Gamma \left( \frac{m}{2}-1 \right) \ \sum_{j=0}^{\infty} \left( 2j + \frac{m}{2} \right) \ \frac{(2j+i)!!}{(2j+m-i-2)!!}\\
& \times  z^{1-m/2} \ J_{2j+m/2}(z) \ C_{2j+1}^{m/2-1}(w)\\
g_m^i(w,z) &=   2^{m/2-1} \ \Gamma \left( \frac{m}{2} \right) \ \sum_{j=0}^{\infty}(4j+m+2) \  \frac{(2j+i)!!}{(2j+m-i)!!}\\
& \times    z^{-m/2} \ J_{2j+m/2+1}(z) \ C_{2j+1}^{m/2}(w).
\end{align*}
\end{theorem}

\begin{proof}
Similar to the even dimensional case, we prove this theorem by induction on the dimension $m$. Hence, we first prove the property for $m=3$. The series expansion for the term $\hat{f}^0_3(s,t) = \widetilde{J}_0(t) = J_0(t)$, namely
\begin{displaymath}
J_0(z \sqrt{1-w^2}) = \sqrt{2} \ \sum_{j=0}^{\infty} \frac{\Gamma \left( j+ \frac{1}{2} \right)}{j!} \ \left( 2j+\frac{1}{2} \right) \ z^{-1/2} \ J_{2j+1/2}(z) \ C_{2j}^{1/2}(w),
\end{displaymath}
can be found in \cite{Erde}, section 7.15, formula (3). The series expansion for $g_3^0$ then follows from the one of $\hat{f}_3^0$ by means of
\begin{displaymath}
g_3^0 = z^{-2}  w^{-1}  \partial_w \hat{f}_3^0 .
\end{displaymath}
Next, the series expansions for $\tilde{f}_3^1$, $\hat{f}_3^1$ and $g_3^1$ follow respectively by use of the formulae
\begin{displaymath}
\tilde{f}_3^1 = - \hat{f}_3^0 , \qquad \hat{f}_3^1 = -s \ \hat{f}_3^0 \qquad \mathrm{and} \qquad g_3^1 = s \ g_3^0  .
\end{displaymath}
The remaining part of the proof is completely similar to the even dimensional case.
\end{proof}
%%%%%%%%%%%%%%%%%%%%%%%%%%%%%%%%%%%%%%%%%%%%%%%%%%%%%%%%%%%%%%%%%%%%%%%%%%%%%%%%%%%%%%%%%%%%%%%%%%%%%%%%%%%%%%%%%%%%%
\subsection{Eigenvalues of new class of Clifford-Fourier transforms}
\label{eigvalsodd}

We will now calculate the action of the new Clifford-Fourier transforms
\begin{displaymath}
\mathcal{F}^{i}_{+,m} \lbrack f(\ux) \rbrack (\uy) = \frac{1}{(2 \pi)^{m/2}} \ \int_{\mathbb{R}^m} K^{i}_{+,m}(\ux,\uy) \ f(\ux) \ dV(\ux), \qquad i=0,1,2,\ldots,m-2
\end{displaymath}
on the basis $\lbrace \psi_{j,k,\ell} \rbrace$. 

From the previous subsection we observe that the new Clifford-Fourier kernels $K_{+,m}^{i}$ in the odd dimensional case are again of the structure (\ref{structure kernel}). Hence, the action of the Clifford-Fourier transforms $\mathcal{F}^{i}_{+,m}$ on the basis $\lbrace \psi_{j,k,\ell} \rbrace$ can once more be determined by substituting the corresponding coefficients $\alpha_k$ and $\beta_k$ in the equations (\ref{eigenvalue eq scaling}). 

\begin{theorem}
\label{EigValsOdd}
In case of $m$ odd, the Clifford-Fourier transforms $\mathcal{F}^{i}_{+,m}$ act as follows on the basis $\lbrace \psi_{j,k,\ell} \rbrace$ of $\mathcal{S}(\mathbb{R}^m) \otimes \cC l_{0,m}$:\vspace{0,2cm}\\
\emph{\underline{Case 1}: $i$ even} ($i=0,2,\ldots,m-5,m-3$).\vspace{0,1cm}\\
\underline{a) $k$: even}
\begin{align*}
\mathcal{F}^{i}_{+,m} \lbrack \psi_{2p,k,\ell} \rbrack (\uy) &=  I \ e_i^c \ (-1)^{(m+1)/2} \ \frac{(k+i-1)!!}{(k-i+m-3)!!} \ (-1)^p \ \psi_{2p,k,\ell}(\uy)\\
\mathcal{F}^{i}_{+,m} \lbrack \psi_{2p+1,k,\ell} \rbrack (\uy) &=  e_i \  \frac{(k+i-1)!!}{(k+m-i-3)!!} \ (-1)^p \ \psi_{2p+1,k,\ell}(\uy)
\end{align*}
\underline{b) $k$: odd}
\begin{align*}
\mathcal{F}^{i}_{+,m} \lbrack \psi_{2p,k,\ell} \rbrack (\uy) &=    e_i \ \frac{(k+i)!!}{(k+m-i-2)!!} \ (-1)^{p+1} \ \psi_{2p,k,\ell}(\uy)\\
\mathcal{F}^{i}_{+,m} \lbrack \psi_{2p+1,k,\ell} \rbrack (\uy) &=  I \ e_i^c \ (-1)^{(m+1)/2} \  \frac{(k+i)!!}{(k-i+m-2)!!} \ (-1)^p \ \psi_{2p+1,k,\ell}(\uy).
\end{align*}
%%%%%%%%%%%%%%%%%%%%%%%%%%%%%%%%%%%%%%%%%%%%%%%%%%%%%%%%%%%%%%%%%%
\emph{\underline{Case 2}: $i$ odd} ($i=1,3,\ldots,m-4,m-2$).\vspace{0,1cm}\\
\underline{a) $k$: even}
\begin{align*}
\mathcal{F}^{i}_{+,m} \lbrack \psi_{2p,k,\ell} \rbrack (\uy) &=   e_i \ \frac{(k+i)!!}{(k+m-i-2)!!} \ (-1)^{p+1} \ \psi_{2p,k,\ell}(\uy)\\
\mathcal{F}^{i}_{+,m} \lbrack \psi_{2p+1,k,\ell} \rbrack (\uy) &=  I \ e_i^c \ (-1)^{(m+1)/2} \  \frac{(k+i)!!}{(k+m-i-2)!!} \ (-1)^{p+1} \ \psi_{2p+1,k,\ell}(\uy)
\end{align*}
\underline{b) $k$: odd}
\begin{align*}
\mathcal{F}^{i}_{+,m} \lbrack \psi_{2p,k,\ell} \rbrack (\uy) &=   I \ e_i^c \ (-1)^{(m+1)/2} \ \frac{(k-1+i)!!}{(k+m-i-3)!!} \ (-1)^{p+1} \ \psi_{2p,k,\ell}(\uy)\\
\mathcal{F}^{i}_{+,m} \lbrack \psi_{2p+1,k,\ell} \rbrack (\uy) &=   e_i \ \frac{(k+i-1)!!}{(k+m-i-3)!!} \ (-1)^p \ \psi_{2p+1,k,\ell}(\uy).
\end{align*}
\end{theorem}

\begin{remark}
Note that for $i=0,1,\ldots,m-3,m-2$ the above eigenvalue equations never reduce to the ones of the Clifford-Fourier transform (see (\ref{eigenvalue C-F})). This means that in the odd-dimensional case, the explicit kernel of the Clifford-Fourier transform is not expressible as a finite sum of powers of $s$ multiplied with Bessel functions in $t$.
\end{remark}
%%%%%%%%%%%%%%%%%%%%%%%%%%%%%%%%%%%%%%%%%%%%%%%%%%%%%%%%%%%%%%%%%%%%%%%%%%%%%%%%%%%%%%%%%%%%%%%%%%%%%%%%%%%%%%%%%
%%%%%%%%%%%%%%%%%%%%%%%%%%%%%%%%%%%%%%%%%%%%%%%%%%%%%%%%%%%%%%%%%%%%%%%%%%%%%%%%%%%%%%%%%%%%%%%%%%%%%%%%%%%%%%%%%
\section{Properties of the new Fourier transforms}
\setcounter{equation}{0}
\label{Properties}

In this section we study some important properties of the integral transforms defined by
\[
\mathcal{F}^{i}_{+,m} \lbrack f(\ux) \rbrack (\uy) = \frac{1}{(2 \pi)^{m/2}} \  \int_{\mathbb{R}^m} K^{i}_{+,m}(\ux,\uy) \ f(\ux) \ dV(\ux)
\]
with kernel $K^{i}_{+,m}$ as given in respectively formula (\ref{kernel even m}) and (\ref{kernel odd m}) for respectively $m$ even and $m$ odd.

We start by obtaining estimates for the kernels.
\begin{lemma}
\label{bounds}
Let $m$ be even and $i= 0, \ldots, m-2$. For $ \ux, \uy \in \RR^m$, there exists a constant $c$ such that
\begin{align*}
|\tilde{f}_m^i(s,t)  + \hat{f}_m^i(s,t) | &\leq  c (1+|\ux|)^{i}(1+|\uy|)^{i},  \\
|(x_{j} y_{k} - x_{k}y_{j}) g_m^i(s,t)| &\leq  c (1+|\ux|)^{i}(1+|\uy|)^{i}, \qquad
 j \neq k. 
\end{align*}
Similarly, in case of $m$ odd and $i=0,\ldots,m-2$, there exists a constant $c$ such that for $ \ux, \uy \in \RR^m$
\begin{align*}
|e_i \ \tilde{f}_m^i(s,t)  + I \ e_i^c \ \hat{f}_m^i(s,t) | &\leq  c (1+|\ux|)^{i}(1+|\uy|)^{i},  \\
|(x_{j} y_{k} - x_{k}y_{j}) \ e_i \  g_m^i(s,t)| &\leq  c (1+|\ux|)^{i}(1+|\uy|)^{i}, \qquad
 j \neq k. 
\end{align*}
\end{lemma}
\begin{proof}
This follows immediately using the well-known bounds
$$
    |  z^{-\a} J_\a(z) | \le c, \qquad z \in  \RR. 
$$
and
$$
       |  z^{-\a + 1} J_\a(z) | \le c, \qquad z \in  \RR, \quad \a \geq \frac{1}{2} 
$$
as in the proof of Lemma 5.2 and Theorem 5.3 in \cite{DBXu}.
\end{proof}

As an immediate consequence of Lemma \ref{bounds}, we can now specify the domain in the definition of the new class of Fourier transforms. Let us define the following function spaces, for $i = 1, \ldots, m-2$,
$$
B_{i}(\RR^m) : = \left\{ f\in L_{1}(\RR^m): \int_{\RR^m} (1+|\uy|)^{i} | f(\uy)| \ dV(\uy) < \infty \right\}.
$$
Note that for $i=0$, $B_{0}(\RR^m) =L_{1}(\RR^m)$.
Then, in the spirit of formulation $\textbf{F1}$ of the ordinary Fourier transform, we have the following theorem.

\begin{theorem} \label{CFdomain}
The integral transform  $\mathcal{F}^{i}_{+,m}$ is well-defined on $B_{i}(\RR^m) \otimes \cC l_{0,m}$. In particular, for $f \in B_{i}(\RR^m) \otimes \cC l_{0,m}$, $\mathcal{F}^{i}_{+,m}\lbrack f \rbrack$ is a continuous function.
\end{theorem}

\begin{proof}
It follows immediately from Lemma \ref{bounds} that the transform is well-defined on $B_{i}(\RR^m) \otimes \cC l_{0,m}$. The continuity of $\mathcal{F}^{i}_{+,m}\lbrack f \rbrack$ follows from the continuity of the kernel and the dominated convergence theorem.
\end{proof}

If we restrict the transforms $\mathcal{F}^{i}_{+,m}$ to the space $\mathcal{S}(\mathbb{R}^m) \otimes \cC l_{0,m}$ of Schwartz class functions taking values in $\cC l_{0,m}$, we can formulate a much stronger result. This is the subject of the following theorem.

\begin{theorem}
\label{TFonS}
Let $i= 0, \ldots, m-2$. The integral transforms $\mathcal{F}^{i}_{+,m}$ define continuous operators mapping $\mathcal{S}(\mathbb{R}^m) \otimes \cC l_{0,m}$ to $\mathcal{S}(\mathbb{R}^m) \otimes \cC l_{0,m}$.

When $m$ is even, the inverse of each transform $\mathcal{F}^{i}_{+,m}$ is given by $\mathcal{F}^{m-2-i}_{+,m}$, i.e.
\[
\mathcal{F}^{i}_{+,m} \mathcal{F}^{m-2-i}_{+,m} = \mathcal{F}^{m-2-i}_{+,m} \mathcal{F}^{i}_{+,m} = id_{\mathcal{S}(\mathbb{R}^m) \otimes \cC l_{0,m}}.
\]
In particular, when $i = (m-2)/2$, the transform reduces to the Clifford-Fourier transform, satisfying
\[
\mathcal{F}^{(m-2)/2}_{+,m} \mathcal{F}^{(m-2)/2}_{+,m} = id_{\mathcal{S}(\mathbb{R}^m) \otimes \cC l_{0,m}}
\]
and the kernel is also given by
\begin{equation}
\label{CFexpr}
K^{(m-2)/2}_{+,m}(\ux,\uy) =  - e^{- \frac{I \pi}{2}\Gamma_{\uy}  } \left( e^{-I \la \ux, \uy \ra} \right).
\end{equation}
\end{theorem}

\begin{proof}
Proving that $\mathcal{F}^{i}_{+,m}$ is a continuous operator on $\mathcal{S}(\mathbb{R}^m) \otimes \cC l_{0,m}$ is done in the same way as in Theorem 6.3 in \cite{DBXu}, so we omit the details.

In case of $m$ even, using the formulas for the eigenvalues (see Theorem \ref{EigValsEven}), we can observe that
\[
\mathcal{F}^{i}_{+,m} \mathcal{F}^{m-2-i}_{+,m} = \mathcal{F}^{m-2-i}_{+,m} \mathcal{F}^{i}_{+,m} = id_{\mathcal{S}(\mathbb{R}^m) \otimes \cC l_{0,m}}
\]
when acting on the eigenfunction basis $\{ \psi_{j,k,l}\}$ of  $\mathcal{S}(\mathbb{R}^m) \otimes \cC l_{0,m}$. As both operators are continuous, the result follows via Hahn-Banach. 

Formula (\ref{CFexpr}) was proven in \cite{DBXu}.
\end{proof}

In the following theorem we discuss the extension of the transforms $\mathcal{F}^{i}_{+,m}$ to $L_{2}(\mathbb{R}^m) \otimes \cC l_{0,m}$.

\begin{theorem}
The transform $\mathcal{F}^{i}_{+,m}$ extends from $\mathcal{S}(\mathbb{R}^m) \otimes \cC l_{0,m}$  to a continuous map on $L_{2}(\mathbb{R}^m) \otimes \cC l_{0,m}$ for all $i  \leq (m-2)/2$, but not for $i > (m-2)/2$. 

In particular, only when $m$ is even and $i = (m-2)/2$, the transform $\mathcal{F}^{(m-2)/2}_{+,m}$ is unitary, i.e.
\[
||\mathcal{F}^{(m-2)/2}_{+,m} (f) || = ||f||
\]
for all $f \in L_{2}(\mathbb{R}^m) \otimes \cC l_{0,m}$.
\end{theorem}

\begin{proof}
The space $L_{2}(\mathbb{R}^m) \otimes \cC l_{0,m}$ is equipped with the inner product 
\[
\langle f, g \rangle = \left[ \int_{\mR^{m}} \overline{f^{c}} \, g \; dV(\ux) \right]_{0}.
\]
Here, the operator $\bar{.}$ is the main anti-involution on the Clifford algebra $\cC l_{0,m}$ defined by
\begin{eqnarray*}
\overline{a b} = \overline{b} \overline{a}, \qquad \overline{e_{i}} = -e_{i}, \quad (i = 1,\ldots, m)
\end{eqnarray*}
and $|.|_{0}$ is the projection on the space of $0$-vectors $\cC l_{0,m}^{0}$. The set of functions $\psi_{j,k,\ell}$ defined in formula (\ref{basis}) is after suitable normalization an orthonormal basis for $L_{2}(\mathbb{R}^m) \otimes \cC l_{0,m}$, satisfying
\[
\langle \psi_{j_{1},k_{1}, \ell_{1}} , \psi_{j_{2},k_{2}, \ell_{2}} \rangle = \delta_{j_{1} j_{2}} \delta_{k_{1} k_{2}}\delta_{\ell_{1} \ell_{2}}, 
\]
see e.g. \cite{MR926831}. 

Now let $f \in L_{2}(\mathbb{R}^m) \otimes \cC l_{0,m}$ have the expansion
\[
f = \sum_{j,k,\ell} a_{j,k,\ell} \psi_{j,k,\ell}
\]
with $\sum_{j,k,\ell} |a_{j,k,\ell}|^{2} < \infty$. Then we compute
\begin{align*}
||\mathcal{F}^{i}_{+,m} (f) ||^{2} &= \sum_{j,k,\ell} |a_{j,k,\ell}|^{2} |\l_{j,k,\ell}|^{2}
\end{align*}
with $\l_{j,k,\ell}$ the eigenvalues of $\mathcal{F}^{i}_{+,m}$ as determined in Theorem \ref{EigValsEven} and \ref{EigValsOdd}.
If $i \leq (m-2)/2$ then $|\l_{j,k,\ell}| \leq 1$ and we have $||\mathcal{F}^{i}_{+,m} (f) || \leq ||f||$ for all $f$. On the other hand, if $i > (m-2)/2$ it is easy to construct an $f \in L_{2}(\mathbb{R}^m) \otimes \cC l_{0,m}$ such that $||\mathcal{F}^{i}_{+,m} (f) || >  \infty$, because then the eigenvalues $\l_{j,k,\ell}$ behave as polynomials in $k$ (when $m$ is even) or as rational functions in $k$ with degree nominator $>$ degree denominator (when $m$ is odd). 

Only when $m$ is even and $i = (m-2)/2$ the eigenvalues have unit norm and the transform is hence unitary.
\end{proof}

We can now also introduce the transforms $\mathcal{F}^{i}_{-,m}$ as
\[
\mathcal{F}^{i}_{-,m} \lbrack f(\ux) \rbrack (\uy) =  \frac{1}{(2 \pi)^{m/2}} \ \int_{\mathbb{R}^m} K^{i}_{-,m}(\ux,\uy) \ f(\ux) \ dV(\ux)
\]
with $K^{i}_{-,m}(\ux,\uy) = \left( K^{i}_{+,m}(\ux,-\uy) \right)^c$. Then we obtain the following proposition
\begin{proposition} \label{DiffTransform}
Let $f \in \cS(\mR^{m})\otimes \cC l_{0,m}$ and $i=0,\ldots,m-2$. Then one has
\begin{align*}
\cF_{\pm,m}^i \left\lbrack \ux \, f \right\rbrack & = \mp \ (\mp I)^m \  \upy \left\lbrack \cF_{\mp,m}^i \lbrack  f \rbrack \right\rbrack\\
\cF_{\pm,m}^i \left\lbrack \upx \lbrack f \rbrack  \right\rbrack & = \mp  \ (\mp I)^m \ \uy \ \cF_{\mp,m}^i \left\lbrack f \right\rbrack.
\end{align*} 
\end{proposition}
\begin{proof}
The first property immediately follows from the first differential equation of the kernel (see (\ref{system}) and (\ref{system2})). The second property follows because of partial integration, which is allowed as the kernel satisfies polynomial bounds (see Lemma \ref{bounds}), and application of the second differential equation of the Clifford-Fourier system.
\end{proof}

\section*{Acknowledgment}

This paper was written when the first author was visiting researcher at the Korteweg-de Vries Institute (University of Amsterdam), supported by a FWO mobility allowance.

%%%%%%%%%%%%%%%%%%%%%%%%%%%%%%%%%%%%%%%%%%%%%%%%%%%%%%%%%%%%%%%%%%%%%%%%%%%%%%%%%%%%%%%%%%%%%%%%%%%%%%%%%%%%%%%%%%%%%%%%
%%%%%%%%%%%%%%%%%%%%%%%%%%%%%%%%%%%%%%%%%%%%%%%%%%%%%%%%%%%%%%%%%%%%%%%%%%%%%%%%%%%%%%%%%%%%%%%%%%%%%%%%%%%%%%

\end{document}